\documentclass[11pt]{amsart}
\usepackage{amssymb, enumerate, mdwlist}

\evensidemargin\paperwidth
\advance\evensidemargin-\textwidth
\advance\oddsidemargin-1in
\evensidemargin\oddsidemargin

\topmargin\paperheight
\advance\topmargin-\textheight
\advance\topmargin-1in

\theoremstyle{plain}
\newtheorem{theorem}{Theorem}[section]
\theoremstyle{plain}
\newtheorem{proposition}[theorem]{Proposition}
\theoremstyle{plain}
\newtheorem{lemma}[theorem]{Lemma}
\theoremstyle{plain}

\theoremstyle{plain}

\theoremstyle{plain}
\newtheorem{mthm}{Theorem}

\theoremstyle{plain}
\newtheorem{mcor}[mthm]{Corollary}
\theoremstyle{definition}
\newtheorem{definition}[theorem]{Definition}
\theoremstyle{remark}
\newtheorem{remark}[theorem]{Remark}
\theoremstyle{remark}
\newtheorem{example}[theorem]{Example}
\theoremstyle{remark}

\title[Group approximation in Cayley topology]
{Group approximation in Cayley topology and coarse geometry, \\ Part I: Coarse embeddings of amenable groups.}
\author{Masato Mimura \and Hiroki Sako}  \thanks{Work was in part done while the first-named author was partially supported  by JSPS KAKENHI Grant Number JP25800033.}
\address{Masato Mimura\\
Mathematical Institute, Tohoku University, Japan\ /\ \'{E}cole Polytechnique F\'{e}d\'{e}rale de Lausanne, Switzerland}
\email{mimura-mas@m.tohoku.ac.jp}
\address{Hiroki Sako\\
School of Science, Department of Mathematical Sciences, Tokai University, Japan}
\email{hiroki.sako@gmail.com}
\date{\today}

\begin{document}

\begin{abstract}
The objective of this series is to study metric geometric properties of (coarse) disjoint unions of amenable Cayley graphs. We employ the Cayley topology and observe connections between large scale structure of metric spaces and group properties of Cayley accumulation points. In this Part I, we prove that a disjoint union has property A of G. Yu if and only if all groups appearing as Cayley accumulation points in the space of marked groups are amenable. 
As an application, 
we construct  two disjoint unions of finite  special linear groups (and unimodular linear groups) with respect to two systems of generators that look similar such that one has property A and the other does not admit (fibred) coarse embeddings into any Banach space with non-trivial type (for instance, any uniformly convex Banach space). 
\end{abstract}


\maketitle

\tableofcontents

\section{Introduction}
The concept of the \textit{coarse embedding} was defined by M. Gromov \cite{GromovAsymptotic}, and can be regarded as an injective morphism in the category of coarse geometry (see Definition~\ref{definition=CoarseEmbedding}; Gromov himself called it a \textit{uniform embedding at infinity}). The concept of \textit{property A} was introduced by G. Yu \cite{Yu}, plays the role of amenability in coarse geometry, and has been regarded as one of the most fundamental properties in coarse geometry (see Definition~\ref{Lemma Property A Hulanicki} in our setting).  Yu proved the following two celebrated properties in \cite{Yu}: For uniformly locally finite (see Definition~\ref{DefinitionLocalFin}) generalized metric spaces $X$ (namely, metric spaces whose metrics possibly take the $+\infty$-value; see Subsection~\ref{subsection=LargeScale}),
\begin{itemize}
  \item if $X$ has property A, then $X$ is coarsely embeddable into a Hilbert space;
  \item if $X$ is coarsely embeddable into a Hilbert space, then the coarse Baum--Connes conjecture holds true for $X$,
\end{itemize}
Here, we make a remark that our convention on coarse embeddability of generalized metric spaces is slightly non-standard. More precisely, we \textit{impose no condition on any pair of points with infinite distance to formulate coarse embeddabilty}; see $(1)$ of Definition~\ref{definition=CoarseEmbedding} for the concrete definition. 
We refer the readers to books \cite{RoeLectureNote} and \cite{BookNowakYu} for comprehensive treatments of these topics.

It has turned out that the class of generalized metric spaces of the form of the \textit{$($coarse$)$ disjoint union of an infinite sequence of finite connected graphs of uniformly bounded degrees} (see Definition~\ref{definition=CoarseDisjointUnion} for disjoint unions and coarse disjoint unions) supplies plenty of examples with notable coarse properties. Here each finite graph is equipped with the shortest path metric. For instance, among uniformly locally finite generalized metric spaces, many examples without property A were given by Guentner (see below); among them, the first example without property A that admits a coarse embeddings into a Hilbert space was provided by Arzhantseva--Guentner--\v{S}pakula \cite{ArzhantsevaGuentnerSpakula}. These examples lie in this class.

The goal of this series of study is to provide a framework to investigate the following class: That of all \textit{$($coarse$)$ disjoint unions of finite connected Cayley graphs} $(\mathrm{Cay}(G_m;S_m))_{m}$ {of uniformly bounded degrees}. There has been extensive study for \textit{box spaces} by various researchers, which are special cases in that class; see below for the precise definition. However, there was no study of the class above in full generality before our result. 

The main idea of us is to employ the concept of \textit{the space of $($$k$-$)$marked groups} $\mathcal{G}(k)$, intensively studied by  R. I. Grigorchuk \cite[Section~6]{Grigorchuk}. For a fixed $k\in \mathbb{N}_{\geq 1}$, $\mathcal{G}(k)$ is the space of all $k$-marked group $\mathbf{G}=(G;s_1,\ldots,s_k)$, where $(s_1,\ldots,s_k)$ is an ordered $k$-tuple that generates $G$. This space is equipped with a natural metrizable and \textit{compact} topology, sometimes called the \textit{Cayley topology} (also known as the Chabauty topology on the space of normal subgroups of $F_k$, the free group of rank $k$; see Remark~\ref{remark=Chabauty}). The convergence in the Cayley topology corresponds to the coincidence of rooted metric balls in \textit{Cayley diagrams} (see Definition~\ref{definition=CayleyDiagram}) of larger and larger radii. We discuss this concept in more details in Subsection~\ref{subsection=SpaceOfMarkedGroups}; see also Lemma~\ref{lemma=Tessera} and the proof of it for an instructive example. For a sequence $(\mathbf{G}_m)_{m\in \mathbb{N}}$ of $k$-marked groups (for a fixed $k$),  we consider two objects of different nature:
\begin{itemize}
  \item The disjoint union $\bigsqcup_{m\in \mathbb{N}}\mathrm{Cay}(\mathbf{G}_m)$ of the Cayley graphs, and the coarse disjoint union $\coprod_{m\in \mathbb{N}}\mathrm{Cay}(\mathbf{G}_m)$: These are generalized metric spaces \textit{without group structure};  
   \item The \textit{Cayley boundary} $\partial_{\mathrm{Cay}}  (\mathbf{G}_m)_m$: It consists of certain $k$-marked \textit{groups}, as we define below.
\end{itemize}
\begin{definition}\label{definition=CayleyBoundary}
For an infinite sequence $(\mathbf{G}_m)_{m\in \mathbb{N}}$ of (possibly infinite) $k$-marked groups, the \textit{Cayley boundary} $\partial_{\mathrm{Cay}} (\mathbf{G}_m)_m$ of it is defined to be the subset of $\mathcal{G}(k)$ consisting of all accumulation points of the sequence $(\mathbf{G}_m)_m$ in $\mathcal{G}(k)$.
\end{definition}
By compactness of $\mathcal{G}(k)$, the Cayley boundary $\partial_{\mathrm{Cay}}  (\mathbf{G}_m)_m$ forms a non-empty \textit{compact} subset of $\mathcal{G}(k)$. 

Our main result(s) build a bridge between \textit{coarse properties} of the former and \textit{group properties} of the latter; see Theorem~\ref{mthm=theoremA} on the topic of the current Part I paper. Before proceeding to that, let us recall what were known for the special case of our class, \textit{box spaces}, including their definition.

\begin{definition}\label{definition=RF}
A finitely generated group $G$ is said to be  \textit{RF} (\textit{Residually Finite}) if there exists a sequence  $(N_m)_{m\in \mathbb{N}}$ of normal subgroups of $G$ of finite index such that $N_{m+1}\leqslant N_m$ and that $\bigcap_{m\in \mathbb{N}} N_m = \{ e_G\}$.
\end{definition}
Fix a finitely generated, infinite, and residually finite group $G$ and fix a finite generating set $S$ of $G$.  Fix a sequence $(N_m)_{m}$ as in Definition~\ref{definition=RF}. Then a coarse disjoint union 
\[
(\Box G =) \Box_{(N_m)_m} G = \coprod_{m \in  \mathbb{N}} \mathrm{Cay}(G/ N_m ;S \ \mathrm{mod}\ N_m)
\]
is called a \textit{box space} of $G$ with respect to  $(N_m)_{m\in \mathbb{N}}$. The coarse structure of $\Box_{(N_m)_m} G$ does heavily depend on the choice of $(N_m)_m$ in general. More precisely, it follows from work of A. Selberg \cite{Selberg} that some box space of $F_2$ forms an expander family \cite{BookLubotzky}. It implies that this box space of $F_2$ does not admit a coarse embedding into a Hilbert space; see also explanation in Proposition~9.7 and Example~9.10 of our Part II paper \cite{MSPartII}. At the other end of the spectrum, as we mentioned above, in \cite{ArzhantsevaGuentnerSpakula} another box space of $F_2$, associated with a different chain $(N_m)_m$ from one in the example of expanders above, is constructed; that box space \textit{does} admit a coarse embedding into a Hilbert space. On the other hand, the choice of $S$ affects metric structures of each component $\mathrm{Cay}(G/ N_m ; S\ \mathrm{mod}\ N_m)$ in biLipschitz ways \textit{only uniformly on $m\in \mathbb{N}$}; therefore, the choice of generating set of $G$ is \textit{irrelevant} in the aspect of coarse properties. E. Guentner  and J. Roe showed the following relation between the coarse geometric properties of $\Box G$ and group properties of the mother group $G$ (see Proposition~11.39 in \cite{RoeLectureNote}); for the definition of \textit{a-}$\mathrm{T}$\textit{-menability}, see Section~\ref{section=a-T-menable}.
\begin{itemize}
\item $\Box G$ has property A \textit{if and only if} $G$ is amenable;
\item $\Box G$ admits a coarse embedding into a Hilbert space \textit{only if} $G$ is a-$\mathrm{T}$-menable. 
\end{itemize}
The former item supplied many examples of (uniformly locally finite) generalized metric spaces without property A; for instance, all box spaces of $F_2$. On the latter assertion, the converse does \textit{not} hold in general; recall a Selberg-type expander family from our discussion above. Chen--Wang--Yu \cite{ChenWangYu} introduced the notion of \textit{fibred coarse embeddings}, which is a certain weak form of (genuine) coarse embeddings, in their study of the maximal coarse Baum--Connes conjecture. Chen--Wang--Wang \cite{ChenWangWang} proved that $\Box G$ admits a fibred coarse embedding into a Hilbert space \textit{if and} only if $G$ is a-$\mathrm{T}$-menable. Fibred coarse embeddings from generalized metric spaces in our class will be studied in our Part II paper \cite{MSPartII}; see Theorem~A and Corollary~B in that paper.

One of our motivations to extend the framework from the class of box spaces to our class is that there exist various natural examples of $((G_m;S_m))_m$ that are \textit{not} derived from box space constructions. Here we exhibit one example.\begin{example}\label{example=SymmetricGroup}
For $m\in \mathbb{N}_{\geq 3}$, let $G_m=\mathrm{Sym}([m])$, the symmetric group over $[m]=\{1,2,\ldots ,m\}$, and $S_m=(\sigma^{(m)},\tau^{(m)})$. Here $\sigma^{(m)}=(12)$ and $\tau^{(m)}=(12\cdots m)$. 
\end{example}
Hereafter, we abbreviate  $\mathrm{Sym}([m])$ as $\mathrm{Sym}(m)$.

To the best knowledge of the authors, it was not known before our result whether the (coarse) disjoint union of $(\mathrm{Cay}(G_m,S_m))_{m\geq 3}$ has property A. We will also consider an analogue of this example for alternating groups; see Remark~\ref{remark=AlternatingGroup}.

In this paper, as main examples, we furthermore deal with several sequences of unimodular linear groups (and special linear groups; see Remark~\ref{remark=SpecialLinearGroup}) $G_m=\mathrm{SL}^{\pm}(m,\mathbb{F}_{p^{n_m}})$ over finite fields with \textit{two different choices of} system of generating sets $(S_m)_{m}$; see Example~\ref{example=SpecialLinear}. Here by \textit{unimodular linear group}, we mean the group of matrices of determinant in $\{\pm 1\}$. We write it as $\mathrm{SL}^{\pm}$. This is an overgroup of index $2$ of the special linear group of the same degree. The relationship between unimodular linear groups and special linear groups resembles that between symmetric groups and alternating groups. We warn that this terminology is \textit{different} from a ``linear group (that means, a group that admits an injective homomorphism into a general linear group of finite degree over a (commutative) field) that is unimodular (in the context of Haar measures on locally compact groups).''

The following is our main theorem in this Part I paper. 
\begin{mthm}
\label{mthm=theoremA}

Let $k\geq 1$ and $(\mathbf{G}_m)_{m\in \mathbb{N}} = \{ ( G_m;s^{(m)}_1,\ldots ,s^{(m)}_k ) \}_{m\in \mathbb{N}}$ be $\mathrm{amenable}$ $($for instance, finite$)$ $k$-marked  groups. 
\begin{enumerate}[$(i)$]
\item 
The disjoint union $\bigsqcup_{m\in \mathbb{N}}\mathrm{Cay} ( \mathbf{G}_m)$  has property A $\mathrm{if\ and\ only\ if}$ every element in $\partial_{\mathrm{Cay}} (\mathbf{G}_m)_m$ is amenable.   
\item 
The disjoint union $\bigsqcup_{m\in \mathbb{N}}\mathrm{Cay} ( \mathbf{G}_m)$ is coarsely embeddable into a Hilbert space $\mathrm{only\ if}$ every element in $\partial_{\mathrm{Cay}} (\mathbf{G}_m)_m$ is a-$\mathrm{T}$-menable.
\end{enumerate}
If $G_m$ is finite for all $m\in \mathbb{N}$, then all of the statements above remain true if we replace the disjoint union $\bigsqcup_{m\in \mathbb{N}}\mathrm{Cay} ( \mathbf{G}_m)$ with a coarse disjoint union $\coprod_{m\in \mathbb{N}}\mathrm{Cay} ( \mathbf{G}_m)$.
\end{mthm}

A rough idea of Theorem~\ref{mthm=theoremA} is that we may regard $\partial_{\mathrm{Cay}} (\mathbf{G}_m)_m$ as the set of ``\textit{all mother groups}'' of the sequence $(\mathbf{G}_m)_m$. Indeed, if $(\mathbf{G}_m)_m=((G/N_m;S\ \mathrm{mod}\ N_m))_m$ comes from a box space, then $\partial_{\mathrm{Cay}} (\mathbf{G}_m)_m$ is the singleton consisting only of the mother group $(G;S)$. 

We emphasize that the following points are visible \textit{only after} extending the framework from the class of box spaces to our general class:
\begin{enumerate}[$(a)$]
  \item The Cayley boundary $\partial_{\mathrm{Cay}} (\mathbf{G}_m)_m$ may consist of \textit{infinitely many} points, possibly it may be a continuum.
  \item Even when $\partial_{\mathrm{Cay}} (\mathbf{G}_m)_m$ is a singleton $\{\mathbf{G}_{\infty}\}$, the Cayley limit group $\mathbf{G}_{\infty}=\lim_{m}\mathbf{G}_m$ is in the class of \textit{LEA}  group when $G_m$, $m\in \mathbb{N}$, is amenable; it is in the class of \textit{LEF} group when $G_m$, $m\in \mathbb{N}$, is furthermore finite. See Definition~\ref{definition=LEF} below for these definitions. In general, the implications
\[
\mathrm{RF}\quad \Longrightarrow \quad \mathrm{LEF}\quad \Longrightarrow \quad \mathrm{LEA} 
\]
hold, and \textit{none} of the implications \textit{can be reversed}. 
  \item Coarse properties of $\bigsqcup_{m\in \mathbb{N}}\mathrm{Cay} ( G_m,S_m)$ may be \textit{considerably affected} by the choice of the system $(S_m)_m$ of generators of $G_m$, even when it might look a slight change; compare with the discussion above on the choice of generating sets $S$ in the case of box spaces.
\end{enumerate}
From these viewpoints, phenomena that the current paper treats provide a complement to recent result of T.\ Delabie and A.\ Khukhro \cite{DelabieKhukhro} and other researchers on coarse geometry of box spaces.

\begin{definition}\label{definition=LEF}
\begin{enumerate}[$(1)$]
 \item (Mal'cev \cite{BookMalcev}, St\"{e}pin \cite{Stepin}, and  Vershik--Gordon \cite{VershikGordon}) A group $G$ is said to be \textit{LEF} (\textit{Locally Embeddable into the class of Finite groups}) if for every finite subset $S\ni e_G$ of $G$, there exists a \textit{finite} group $H$ and a \textit{partial homomorphism} $\beta\colon S\to H$ that is injective. Here, $\beta$ is called a \textit{partial homomorphism} if the following holds true:
\[
\textrm{For every $g,h\in S$ } \textit{such that $gh\in S$}\textrm{,\quad $\beta(gh)=\beta(g)\beta(h)$.}
\]

For a finitely generated group $G$, this definition is equivalent to saying that $\mathbf{G}$ is in the closure (in the Cayley topology) of all finite marked groups for some (equivalently, all) marking $\mathbf{G}$ of $G$.
 \item (Gromov \cite{GromovLEA}) A group $G$ is said to be \textit{LEA} (\textit{Locally Embeddable into the class of Amenable groups}) if for every finite subset $S\ni e_G$ of $G$, there exists an \textit{amenable} group $H$ and an injective partial homomorphism $\beta\colon S\to H$. For a finitely generated group $G$, this definition is equivalent to saying that $\mathbf{G}$ is in the closure (in the Cayley topology) of all amenable marked groups for some (equivalently, all) marking $\mathbf{G}$ of $G$.
\end{enumerate}
In $(2)$, Gromov himself used the terminology of the  \textit{initial subamenability}.
\end{definition}

For concrete examples that illustrate points $(a)$ and $(b)$, see, respectively, Example~\ref{example=Grigorchuk} and Example~\ref{example=Point(b)}. Here we exhibit examples on point $(c)$.

\begin{example}\label{example=SpecialLinear}
Fix  a prime $p$. For $n\in \mathbb{N}_{\geq 1}$, denote by $\mathbb{F}_{p^n}$ the finite field of order $p^n$. It is well known that the multiplicative group $\mathbb{F}_{p^n}^{\times}$ is cyclic; for each $p$ and each $n$, we fix a generator $t_n=t_{p,n}\in \mathbb{F}_{p^n}$  of it. Fix a sequence $(n_m)_{m\in \mathbb{N}_{\geq 3}}$ of positive integers such that $\lim_{m\to \infty}n_m=+\infty$.

Let $G_m$ be the unimodular linear group over $\mathbb{F}_{p^{n_m}}$ of degree $m$, namely,
\[
G_m=\mathrm{SL}^{\pm}(m,\mathbb{F}_{p^{n_m}})=\{g\in \mathrm{GL}(m,\mathbb{F}_{p^{n_m}}):\mathrm{det}(g)\in \{\pm 1\}\}.
\]
Then we consider the following two systems $(S_m)_{m\in \mathbb{N}_{\geq 3}}$, $(T_m)_{m\in \mathbb{N}_{\geq 3}}$ of generators of $G_m$; see Remark~\ref{remark=Generators} for the proof of the fact that $S_m$ and $T_m$ both generate $G_m$.

\begin{itemize}
  \item For $m\in \mathbb{N}_{\geq 3}$, $S_m=(\sigma^{(m)},\upsilon^{(m)},\delta^{(m)},\tau^{(m)})$. Here 
\[
\sigma^{(m)} = 
\left(
\begin{array}{ccccc}
1 & 1 & 0 & \cdots & 0 \\
0 & 1 & 0 &    \        & \vdots \\
\vdots & 0 & 1 & \ddots & \vdots \\
\vdots & \ & \ddots & \ddots & 0 \\
0 & \cdots & \cdots & 0 & 1 \\
\end{array}
\right),
\quad 
\upsilon^{(m)} = 
\left(
\begin{array}{ccccc}
1 & t_{n_m} & 0 & \cdots & 0 \\
0 & 1 & 0 &     \       & \vdots \\
\vdots & 0 & 1 & \ddots & \vdots \\
\vdots & \ & \ddots & \ddots & 0 \\
0 & \cdots & \cdots & 0 & 1 \\
\end{array}
\right),
\]
\[
\delta^{(m)}=
\left(
\begin{array}{ccccc}
-1 & 0 & \cdots & \cdots & 0 \\
0 & 1 & 0 &      \      & \vdots \\
\vdots & 0 & 1 & \ddots & \vdots \\
\vdots & \ & \ddots & \ddots & 0 \\
0 & \cdots & \cdots & 0 & 1 \\
\end{array}
\right), \quad \textrm{and} \quad
\tau^{(m)} = 
\left(
\begin{array}{ccccc}
0 & \cdots & \cdots & 0& 1 \\
1 & 0         &  \          &  \  & 0 \\
0 & 1         & 0         &  \  & \vdots \\
\vdots & \ddots & \ddots & \ddots & \vdots \\
0 & \cdots & 0 & 1 & 0 \\
\end{array}
\right).
\]
  \item For $m\in \mathbb{N}_{\geq 3}$, $T_m=(\sigma^{(m)},\sigma'^{(m)},\upsilon^{(m)}, \delta^{(m)},\tau^{(m)})$. Here $\sigma^{(m)}$, $\upsilon^{(m)}$, $\delta^{(m)}$, and $\tau^{(m)}$ are the same as above, and $\sigma'^{(m)}={}^t \sigma^{(m)}$ is the transpose of $\sigma^{(m)}$.
\end{itemize}
\end{example}

Then, Theorem~\ref{mthm=theoremA} yields the following corollary.

\begin{mcor}\label{mcor=ChoiceOfGenerators}
\begin{enumerate}[$(i)$]
  \item In Example~$\ref{example=SymmetricGroup}$, the disjoint union 
\[
X=\bigsqcup_{m\in \mathbb{N}_{\geq 3}} (\mathrm{Sym}(m);\sigma^{(m)},\tau^{(m)})
\]
has property A.
 \item In Example~$\ref{example=SpecialLinear}$, the disjoint union 
\[
Y=Y_{p,(n_m)_m}=\bigsqcup_{m\in \mathbb{N}_{\geq 3}} (\mathrm{SL}^{\pm}(m,\mathbb{F}_{p^{n_m}});\sigma^{(m)},\upsilon^{(m)},\delta^{(m)},\tau^{(m)})
\]
has property A.
 \item In Example~$\ref{example=SpecialLinear}$, the disjoint union 
\[
Z=Z_{p,(n_m)_m}=\bigsqcup_{m\in \mathbb{N}_{\geq 3}} (\mathrm{SL}^{\pm}(m,\mathbb{F}_{p^{n_m}});\sigma^{(m)},\sigma'^{(m)},\upsilon^{(m)},\delta^{(m)},\tau^{(m)})
\]
does $\mathrm{not}$ admit coarse embeddings into a Hilbert space.
\end{enumerate}
The same assertions, respectively, remain true if we replace the disjoint unions with coarse disjoint unions.
\end{mcor}

\begin{remark}\label{remark=NontrivialType}
In fact, in $(iii)$ above, $Z$ does \textit{not} admit a coarse embedding into any Banach space with non-trivial (linear, also known as Rademacher) type; the class of all Banach spaces of non-trivial type contains all \textit{uniformly convex} Banach spaces. See \cite{BookTomczak-Jaegermann} and \cite{BookBenyaminiLindenstrauss} for the definitions of these concepts; see also Example~4.11 in our Part II paper \cite{MSPartII}. This non-embeddability result follows from a combination of results of V. Lafforgue \cite{Lafforgue1}, \cite{Lafforgue2} on strong Banach property $(\mathrm{T})$ and of arguments in our Part II paper; see \cite[Corollary~B and Corollary~2.2]{MSPartII}.
\end{remark}

In study of constructing \textit{expanders} (see \cite{BookLubotzky} for this topic) from a sequence of finite groups, the phenomenon as in $(c)$ (heavy dependence of the choice of the system of generators) has been observed by various researchers. Nevertheless, we expect that our examples in Corollary~\ref{mcor=ChoiceOfGenerators} might be of its own interest; we only add $\sigma'^{(m)}$ to $S_m$ to construct $T_m$. It might seem as if this perturbation were small, but \textit{it}, in fact, \textit{completely destroys flexibility} (such as property A) of $Y$. To observe this destruction, we might need to come up with the notion of the Cayley limit groups of $((G_m,S_m))_m$ and $((G_m,T_m))_m$. 

\begin{remark}\label{remark=OtherRings}
The choices of coefficient rings of unimodular linear groups in Example~\ref{example=SpecialLinear} does \textit{not} play the most important role in Corollary~\ref{mcor=ChoiceOfGenerators} (only the changes of ``limit (marked) ring'' are slightly involved). The most important points in the construction are the \textit{absence} of $\sigma'^{(m)}$ in $S_m$ and the \textit{presence} of that in $T_m$; see the proof of Proposition~\ref{proposition=LimitLinear} and the discussion below  the proof of Proposition~$\ref{proposition=ChoiceOfGenerators}$. In fact, we construct another example with coefficient rings being $\mathbb{Z}/l_m\mathbb{Z}$, where $l_m\to \infty$, in a very similar manner; this may be regarded as an archimedean analogue of Example~\ref{example=SpecialLinear}. See Example~\ref{example=Point(c)} for that example and Proposition~\ref{proposition=ChoiceOfGenerators} for the conclusion on its coarse properties. See discussion in the last paragraph of Section~\ref{section=PoorCompression} for yet one more construction.

We also construct examples for special linear groups; see Remark~\ref{remark=SpecialLinearGroup}.
\end{remark}

Finally, we prove that a disjoint union $Y$ in item $(ii)$ in Corollary~\ref{mcor=ChoiceOfGenerators} has very poor metric \textit{compression function} (see Definition~\ref{definition=CoarseEmbedding}) for a suitable choice of $(n_m)_{m\geq 3}$. Recall from Theorem~\ref{mthm=theoremA} that $Y$ there has property A; from the Dvoretzky theorem \cite[Chapter~12]{BookBenyaminiLindenstrauss}, this implies that for every infinite dimensional Banach space $E$, $Y$ admits a coarse embedding into $E$. 

\begin{mthm}
\label{mthm=PoorCompression}
In let $Y=Y_{p,(n_m)_m}$ be as in $(ii)$ of Corollary~\ref{mcor=ChoiceOfGenerators}. 

Then, for every non-decreasing proper function $\rho \colon [0,\infty) \to [0,\infty)$, there exists $($explicit$)$ $(n_m)_m$ such that for every prime $p$, $Y=Y_{p,(n_m)_m}$ has the following property: ``Let $E$ be a complex Banach space that is sphere equivalent $($see below$)$ to a complex Banach space $F$ of non-trivial $($linear$)$ type. Then, $(\rho(r),\omega(r)=r)$ is $\mathrm{not}$ a control pair $($see Definition~$\ref{definition=CoarseEmbedding}$$)$ for $Y_{p,(n_m)_m}$ into $E$.''

In particular, it holds for all complex uniformly convex Banach spaces $E$.
\end{mthm}

Here, two Banach spaces $E$ and $F$ are said to be \textit{sphere equivalent} if there exists a bijection $\Phi$ between the unit spheres $S(E)$ and $S(F)$ such that $\Phi$ and $\Phi^{-1}$ are both uniformly continuous (see Section~\ref{section=PoorCompression}). One remarkable example is that $L_{1,\mathbb{C}}=L_1([0,1],\mathbb{C})$ is sphere equivalent to $L_2$; the latter is uniformly convex but the former is of trivial type (see \cite[Chapter~9]{BookBenyaminiLindenstrauss}). Similar holds even in the context of non-commutative $L_p$ spaces; see \cite{Ricard}. It is well known that $L_{1,\mathbb{C}}$ admits a coarse embedding into $L_2$. However, the analogue is \textit{no longer true} in the setting of non-commutative $L_p$ spaces; see Corollary~8.17 and discussions below in \cite{BookBenyaminiLindenstrauss}. Note that a coarse embedding from the disjoint union of uniformly discrete coarse geodesic  metric spaces  must be Lipschitz. Since every Banach space is isometrically isomorphic to its rescaling, without loss of generality we may assume that  $\omega(r)=r$, as we did in Theorem~\ref{mthm=PoorCompression}.

After the first version of our draft came out, remarkably, V. Alekseev and M. Finn-Sell \cite{AlekseevFinnSell} extended the framework of LEF approximation by finite Cayley diagrams to sofic approximations. Nevertheless, we point out that our original setting might remain to be of importance for the following reasons. First, sofic approximation does not cover the case of LEA approximation by infinite amenable Cayley diagrams. Secondly, the framework in the present paper supplies concrete examples of generalized metric spaces with notable features as we saw in Corollary~\ref{mcor=ChoiceOfGenerators} and Theorem~\ref{mthm=PoorCompression}; see also Proposition~\ref{proposition=ChoiceOfGenerators} and examples in Subsection~\ref{subsection=Examples} and Remarks~\ref{remark=SymmetricGroup} and \ref{remark=LinearGroup}. Thirdly, on (fibred) coarse emdeddings into metric spaces that are not Hilbert spaces, even though our theory applies, there would be difficulty in the setting of sofic approximations. See also Point $(a)$ above. 

We also note that in the general framework of \cite{AlekseevFinnSell}, one direction of $(i)$ of Theorem~\ref{mthm=theoremA} (property A of the approximation imples amenability of the limit group) was deduced; see work of T. Kaiser \cite{Kaiser} for a counterexample to the converse direction and for further developments.


\ 

\paragraph{\bf Notation and Conventions:} We use $G$ for a (non-marked) group and $\mathbf{G}$ for a marked group. We write the group unit of $G$ as $e_G$. A finite generating set $S$ of $G$ is regarded as an ordered set (sometimes an ordered multi-set) $S=(s_1,s_2,\ldots ,s_k)$ so that $(G;S)$ is seen as a marked group. A marked group $\mathbf{G}=(G;S)$ is said to be finite (respectively, amenable, and  a-$\mathrm{T}$-menable) if so is $G$. For $k\in \mathbb{N}_{\geq 1}$, we denote by $\mathbf{F_k}$ the \textit{free $k$-marked group}, namely, $\mathbf{F_k}=(F_k;a_1,\ldots ,a_k)$. Here $(a_1,\ldots ,a_k)$ denotes a free basis of $F_k$. For a non-empty set $B$, denote by $\mathrm{Sym}(B)$ the full symmetric group, and by $\mathrm{Sym}_{<\aleph_0}(B)$ the symmetric group with finite support, namely, the group of all permutations on $B$ that fix all but finitely many elements in $B$. For $m\in \mathbb{N}_{\geq 1}$, let $[m]=\{1,2,\ldots,m\}$. For $m,n \in \mathbb{Z}$ with $n\leq m$, set  $[n.m]=\{n,n+1,\ldots,m-1,m\}$. For a set $B$, denote by $\ell_2(B)$ the (Hilbert) space of $\ell_2$-functions $B\to \mathbb{C}$.

In Section~\ref{section=PoorCompression}, we employ some symbols to indicate orders. For two non-negative valued function $a$ and $b$ defined on the same parameter set $\mathcal{P}$, we write as $a \precsim b$ (or $b\succsim a$) if there exists a \textit{universal} constant $C > 0$  such that for all $p \in \mathcal{P}$, $a(p) \leq C b(p)$ holds true. We use the symbol $a \asymp b$ if both $a \precsim b$ and $a \succsim b$  hold. For variables $p_1,\ldots ,p_l$ in  $\mathcal{P}$, we write as $a \precsim_{p_1,\ldots ,p_l} b$  if the positive multiplicative constant $C=C(p_1,\ldots ,p_l)$ above (may not be universal but) may depend on $p_1,\ldots ,p_l$ but it does not depend on the other parameters in $\mathcal{P}$. We write $a \precnsim b$ if $a\precsim b$ holds but $a \succsim b$ fails to be true.

\

\paragraph{\bf Organization of the paper:} In Section~\ref{section=CayleyTopology}, we describe the definition and basic properties of the space of marked groups and the Cayley topology; in Subsection~\ref{subsection=Examples}, we exhibit various examples that illustrate points $(a)$, $(b)$ and $(c)$ above. Section~\ref{section=Amenability} is devoted to the proof of $(i)$ of Theorem~\ref{mthm=theoremA}. In section~\ref{section=a-T-menable}, we prove  $(ii)$ of Theorem~\ref{mthm=theoremA}. In Section~\ref{section=Examples}, we identify the Cayley limit groups in concrete situations and prove Corollary~\ref{mcor=ChoiceOfGenerators} and Proposition~\ref{proposition=ChoiceOfGenerators}. In Section~\ref{section=PoorCompression}, we prove Theorem~\ref{mthm=PoorCompression}.


\section{Preliminaries}\label{section=CayleyTopology}
\subsection{The space of $k$-marked groups and the Cayley topology}\label{subsection=SpaceOfMarkedGroups}
We briefly describe the definition of $\mathcal{G}(k)$ and the Cayley topology. For more detail, see \cite[Section~6]{Grigorchuk} and \cite[Section V.10]{BookdelaHarpe}.

Fix $k\in \mathbb{N}_{\geq 1}$. A $k$\textit{-marked group} $\mathbf{G}=(G;S)=(G;s_1,\ldots,s_k)$ is the pair of a group $G$ and an ordered $k$-tuple $S=(s_1,\ldots , s_k)$ of generators (as a group). We identify two marked groups $(G; s_1, s_2, \cdots, s_k)$ and $(G^\prime; s_1^\prime$, $s_2^\prime, \cdots, s_k^\prime)$ if there exists an isomorphism $\varphi \colon G \stackrel{\simeq}{\rightarrow} G^\prime$ satisfying $\varphi(s_j) = s_j^\prime$ for all $j\in [k]$. 
We denote by $\mathcal{G}(k)$ the set of all (the isomorphism classes of) marked groups.

The key observation is that there is a natural one to one correspondence
\[
\mathcal{G}(k) \quad \longleftrightarrow \quad \{N:N\trianglelefteq F_k\} \quad (\subseteq \{0,1\}^{F_k})
\]
by $\mathbf{G}=(G;s_1,\ldots ,s_k) \mapsto \mathrm{Ker}\{\varphi_{\mathbf{G}}\colon\mathbf{F_k}\twoheadrightarrow \mathbf{G}\}$. Here, $\varphi_{\mathbf{G}}$ above is given by sending $a_j$ to $s_j$ for every $j\in [k]$. (Recall our notation of the free $k$-marked group $\mathbf{F_k}=(F_k;a_1,\ldots,a_k)$ from Introduction.)

\begin{definition}\label{definition=CayleyTopology}
The \textit{Cayley topology} on $\mathcal{G}(k)$ is defined as the relative topology on the set $\{N:N\trianglelefteq F_k\}$, via the identification above, induced by the product topology on $\{0,1\}^{F_k}$.

We write the convergence of a sequence in this topology as $\mathbf{G}_m\stackrel{\mathrm{Cay}}{\to}\mathbf{G}_{\infty}$ or as $\lim_{m\to \infty}\mathbf{G}_m=\mathbf{G}_{\infty}$.
\end{definition}
Since the set $\{N:N\trianglelefteq F_k\}$ is closed in $\{0,1\}^{F_k}$, the Cayley topology is metrizable and \textit{compact}.

There are two viewpoints of the Cayley topology: One is \textit{logical} (relations on groups) and the other one is \textit{geometric} (balls of Cayley diagrams). Before proceeding to them, we give the definition of \textit{Cayley diagrams} and \textit{Cayley graphs} of marked groups; \textit{these two notions are distinguished in this paper.}

\begin{definition}\label{definition=CayleyDiagram}
Let $\mathbf{G}=(G;S)=(G;s_1,\ldots,s_k) \in \mathcal{G}(k)$.
\begin{enumerate}[$(i)$]
\item \begin{enumerate}[$(1)$]
  \item The \textit{Cayley diagram}, written as $\mathrm{CayD}(\mathbf{G})$, is an \textit{edge-colored} and \textit{edge-oriented} graph constructed as follows: The vertex set is $G$. Each edge is colored in either of $k$ colors $[k]=\{1,2,\ldots ,k\}$. For every $g\in G$ and $s_j$, $j\in [k]$, draw an edge from $g$ to $s_jg$;  we color this edge in color $j(\in [k])$ and put an orientation from $g$ to $s_jg$. The edge set is constructed by performing this procedure for all $g\in G$ and $s_j$ for all $j\in [k]$. The resulting $\mathrm{CayD}(\mathbf{G})$  is a diagram, possibly with self-loops or multiple edges.
  \item The \textit{Cayley graph}, written as $\mathrm{Cay}(\mathbf{G})$, is a graph obtained by forgetting both of edge-colorings and edge-orientations from $\mathrm{CayD}(\mathbf{G})$.
\end{enumerate}
Both of $\mathrm{CayD}(\mathbf{G})$ and $\mathrm{Cay}(\mathbf{G})$ are equipped with the shortest path metric on vertices, denoted by $d_{\mathbf{G}}$; for $\mathrm{CayD}(\mathbf{G})$, we ignore edge-orientations for paths. Hence, $\mathrm{CayD}(\mathbf{G})$ and $\mathrm{Cay}(\mathbf{G})$ are isometric as metric spaces.
\item \begin{enumerate}[$(1)$]
  \item For $R\in \mathbb{N}$ and for a non-empty subset $Y$ of $G$ let ${B}_{\mathbf{G}}(Y, R)$ denote the neighborhood 
\[
{B}_{\mathbf{G}}(Y, R) = \{g \in G : \textrm{there exists $h \in Y$ such that $d_{\mathbf{G}}(g, h) \le R$}\}
\]
and let $\partial_{\mathbf{G}}(Y,R)$ denote the boundary ${B}_{\mathbf{G}}(Y, R) \setminus Y$. For a singleton $\{g\}$, we simply write ${B}_{\mathbf{G}}(g, R)$ for ${B}_{\mathbf{G}}(\{g\}, R)$.
  \item For $R\in \mathbb{N}$ and for $g\in G$, let $B_{\mathrm{CayD}(\mathbf{G})}(g,R)$ denote the closed ball of radius $R$ centered at $g$ \textit{as a diagram}, namely, we take the induced graph of $\mathrm{CayD}(\mathbf{G})$ on the vertex set ${B}_{\mathbf{G}}(g, R)$ and construct a diagram by keeping all information on edge-colorings and edge-orientations from the original diagram. We regard $g$ as the root of $B_{\mathrm{CayD}(\mathbf{G})}(g,R)$; thus, $B_{\mathrm{CayD}(\mathbf{G})}(g,R)$ is in fact a \textit{rooted} diagram.
\end{enumerate}
\end{enumerate}
\end{definition}

By construction, each vertex $g$ of $\mathrm{CayD}(\mathbf{G})$ has degree $2k$: One edge outgoing from $g$ in color $j$ for each $j\in [k]$, and one edge incoming to $g$ in color $j$ for each $j\in [k]$. Here if there exists a self-loop (this amounts to saying that $s_j=e_G$ for some $j\in [k]$), then we consider this loop as two edges: One outgoing and the other incoming. Consider two diagrams with the same edge-coloring set $[k]$. We then say a map between them is a \textit{diagram automoprhism} if it is a graph isomorphism and if it furthermore \textit{preserves edge-colorings and edge-orientations}.

There are two ways to define Cayley diagrams/graphs. In this paper, we adopt that of \textit{right} Cayley diagrams/graphs. Namely, $G$ acts on $\mathrm{CayD}(\mathbf{G})$ as diagram automoprhisms transitively \textit{from the right}, and it acts on $\mathrm{Cay}(\mathbf{G})$ as graph automorphisms from the right. The metric on $\mathrm{CayD}(\mathbf{G})$ (on $\mathrm{Cay}(\mathbf{G})$) coincides with the \textit{right-invariant} word metric on $G$ with respect to the marking $S$; as we did in Definition~\ref{definition=CayleyDiagram}, we write that metric as $d_{\mathbf{G}}$ in this paper. In this way, we view $\mathbf{G}$ as a geometric object. 

\begin{example}\label{example=NotIsomorphicAsDiagram}
The Cayley diagram $\mathrm{CayD}(\mathbf{G})$ keeps complete information of $\mathbf{G}$. On the other hand, the Cayley graph $\mathrm{Cay}(\mathbf{G})$ loses much information in general. For instance, let $G$ be a finite group and $S_G=G\setminus \{e_G\}$ (with some ordering). Then, for $\mathbf{G}=(G;S_G)$, the Cayley graph $\mathrm{Cay}(G;S_G)$ is the graph obtained by doubling each edge of the complete graph $K_{\sharp (G)}$ of size $\sharp (G)$. Hence, whenever $\sharp (G)=\sharp (H)<\infty$ (even when $G$ is simple and $H$ is abelian), $\mathrm{Cay}(G;S_G)$ and $\mathrm{Cay}(H;S_H)$ are isomorphic (as graphs). 

To get an example among infinite groups, for instance, take wreath products with a fixed infinite (marked group); see Example~\ref{example=Point(b)} for the definition.
\end{example}

We now explain two viewpoints of the Cayley topology.

\begin{enumerate}[$(1)$]
  \item (\textit{Logical} viewpoint)
        A \textit{relation} $r$ is an element in $F_k$. For a relation $r$, each marked group $\mathbf{G}$ is classified into either of two cases: One is that the aforementioned surjection $\varphi_{\mathbf{G}}\colon \mathbf{F_k}\twoheadrightarrow \mathbf{G}$ sends $r$ to $e_G$ (in this case, we say that $\mathbf{G}$ \textit{satisfies} the relation $r$), the other is that $\varphi_{\mathbf{G}}$ sends $r$ to a non-neutral element (in that case, we say that $\mathbf{G}$ \textit{fails} $r$). In this way, $\mathcal{G}(k)$ is decomposed into 
\[
\mathcal{G}(k)=\{\mathbf{G}:\textrm{$\mathbf{G}$ satisfies $r$}\}\sqcup\{\mathbf{G}:\textrm{$\mathbf{G}$ fails $r$}\}.
\]
These two sets in the right hand side are both \textit{clopen} (closed and open) in the Cayley topology. Moreover, these sets for all $r\in F_k$ generate the topology.
  \item (\textit{Geometric} viewpoint)
    The convergence in the Cayley topology may be regarded as the \textit{local convergence} (equivalently, the \textit{Gromov--Hausdorff convergence} in the current setting) \textit{among rooted diagrams}. Precise meaning is as follow.

\begin{align*}
&\textrm{$(\mathbf{G}_m)_{m\in \mathbb{N}}$ converges to $\mathbf{G}_{\infty}$ in $\mathcal{G}(k)$ (in the Cayley topology)}\tag{$\star$}\\
\Longleftrightarrow \quad &\textrm{``For every $R\in \mathbb{N}$, there exists $m_R\in \mathbb{N}$ such that for all $m\geq m_R$,}\\
&\textrm{ $B_{\mathrm{CayD}(\mathbf{G}_m)}(e_{G_m},R)\cong B_{\mathrm{CayD}(\mathbf{G}_{\infty})}(e_{G_{\infty}},R)$ \textit{as rooted diagrams.''}} 
\end{align*}
Here an isomorphism $B_{\mathrm{CayD}(\mathbf{G}_m)}(e_{G_m},R)\stackrel{\cong}{\to} B_{\mathrm{CayD}(\mathbf{G}_{\infty})}(e_{G_{\infty}},R)$ \textit{as rooted diagrams} means a diagram isomorphism that preserves roots ($e_{G_m}\mapsto e_{G_{\infty}}$). From our discussion below Definition~\ref{definition=CayleyDiagram}, such an isomorphism is canonical; it is unique if exists.
\end{enumerate}

In the viewpoint of $(1)$, $(\star)$ amounts to saying the following. 
\begin{lemma}\label{lemma=Neighborhood}
Let $\mathbf{G}_m=(G_m;s_1^{(m)},\ldots ,s_k^{(m)})$, $m\in \mathbb{N}$, and $\mathbf{G}_{\infty}=(G_{\infty};s_1^{(\infty)},\ldots ,s_k^{(\infty)})$ be in $\mathcal{G}(k)$. Then, $\lim_{m\to \infty}\mathbf{G}_m=\mathbf{G}_{\infty}$ if and only if for every $R\in \mathbb{N}$, there exists $m_R\in \mathbb{N}$ such that for all $m\geq m_R$, the map defined by sending $s_j^{(m)}$ to $s_j^{(\infty)}$ for every $j\in [k]$ gives rise to a $\mathrm{partial}$ $\mathrm{isomorphism}$ $\beta_{\mathbf{G}_m,\mathbf{G},R}\colon B_{\mathbf{G}_m}(e_{G_m},R)\to B_{\mathbf{G}_{\infty}}(e_{G_{\infty}},R)$. Here a partial isomorphism means a  bijective partial homomorphism $($recall Definition~$\ref{definition=LEF}$$)$.

In other words, for $\mathbf{G}=(G;s_1,\ldots ,s_k)\in \mathcal{G}(k)$, if we define for each $R\in \mathbb{N}$,
\begin{align*}
N(\mathbf{G},R)=\{&\mathbf{H}=(H;t_1,\ldots ,t_k)\in \mathcal{G}(k): 
\textrm{the map $t_j\mapsto s_j$ induces}\\ 
&\textrm{a partial isomorphism} \quad \beta_{\mathbf{H},\mathbf{G},R}\colon B_{\mathbf{H}}(e_{H},R) \to B_{\mathbf{G}}(e_{{G}},R).\},
\end{align*}
then $\{N(\mathbf{G},R)\}_{R\in \mathbb{N}}$ forms an $($open$)$ neighborhood system of $\mathbf{G}$.
\end{lemma}

We describe some intuition of the \textit{local} convergence by a pedagogical example in Subsection~\ref{subsection=Local}. We suggest the readers consult Lemma~\ref{lemma=Tessera} and the proof of it.

\begin{remark}\label{remark=Chabauty}
The Cayley topology may be seen as (some restriction of) the \textit{Chabauty topology} (for a locally compact group) in a special case. More precisely, the Chabauty topology for $F_k$ is the topology on the set of all subgroups of $F_k$, which is defined as the relative topology to the product topology on $\{0,1\}^{F_k}$; the Cayley topology on $\mathcal{G}(k)$ is identical to the Chabauty topology above, restricted to the set of all \textit{normal} subgroups of $F_k$.
\end{remark}

For two marked groups $\mathbf{G}=(G;s_1,\ldots,s_k)$ and $\mathbf{G}^{\prime}=(G^{\prime};s_1^{\prime},\ldots,s_k^{\prime})$, we say that $\mathbf{G}'$ is a \textit{marked group quotient} of $\mathbf{G}$ if there exists a surjective homomorphism $G\twoheadrightarrow G^{\prime}$ sending $s_j$ to $s_j^{\prime}$ for every $j\in [k]$. We say that $\mathbf{G}$ is \textit{finitely presented} the corresponding kernel of $\mathbf{F_k}\twoheadrightarrow \mathbf{G}$ is generated as a normal subgroup of $F_k$ by finitely many elements. The finite presentability of $\mathbf{G}$ does not depend on the choice of markings; see \cite[V.2]{BookdelaHarpe}. Recall the definitions of the Cayley boundary and the LEF property, respectively, from Definitions~\ref{definition=CayleyBoundary} and \ref{definition=LEF}. The following lemma was proved by Vershik--Gordon \cite[2.~Theorem]{VershikGordon}.

\begin{lemma}\label{lemma=LEFandRF}
\begin{enumerate}[$(1)$]
  \item All finite marked groups are isolated in the space of marked groups. In particular, if $(\mathbf{G}_m)_{m\in \mathbb{N}}$ is a sequence of finite marked groups in $\mathcal{G}(k)$ such that $\lim_{m\to \infty} \sharp(G_m)=\infty$, then 
\[
\partial_{\mathrm{Cay}}(\mathbf{G}_m)_m=\overline{\{\mathbf{G}_m:m \in \mathbb{N}\}}^{\mathrm{Cay}} \setminus \{\mathbf{G}_m:m \in \mathbb{N}\}.
\]
Here $\overline{\{\ \}}^{\mathrm{Cay}}$ denotes the closure in the Cayley topology.
  \item Let $\mathbf{G}$ be a finitely presented marked group. Then, 
\[
\mathcal{Q}_{\mathbf{G}}=\{\textrm{marked group quotients of $\mathbf{G}$}\}
\]
forms a $($closed and$)$ $\mathrm{open}$ set in the Cayley topology. In particular, for a $\mathrm{finitely}$ $\mathrm{presented}$ group $G$, $G$ is RF if and only if it is LEF.
\end{enumerate}
\end{lemma}

Before proceeding to the proof of Lemma~\ref{lemma=LEFandRF}, we make the following remark. In the definition of box spaces (see below Definition~\ref{definition=RF}), we assume that $(N_m)_m$ is nested (in other words, $N_{m+1}\leqslant N_m $ holds for every $m\in \mathbb{N}$). In general, we may remove this condition: for general sequence $(N_m)_{m\in \mathbb{N}}$ of finite index normal subgroups, if we set $N'_m=\bigcap_{l\leq m} N_l$, then $(N'_m)_{m\in \mathbb{N}}$ is a \textit{nested} sequence of finite index normal subgroups. However, if we drop the condition of being nested, we need to modify the other condition $\bigcap_{m\in \mathbb{N}}N_m=\{e_G\}$ to
\[
\mathrm{liminf}_{m\to \infty}N_m(=\bigcup_{m\in \mathbb{N}}\bigcap_{n\in \mathbb{N}_{\geq m}}N_n)=\{e_G\}.
\]
These two conditions are equivalent if $(N_m)_m$ is nested. However, in general case of a sequence of finite index normal subgroups, the latter condition is the right one to ensure that $((G/N_m;S\ \mathrm{mod}\ N_m))_m$ converges to $(G;S)$ in the space of marked groups. Compare with the discussion in Remark~\ref{remark=FromLEFToRF}.

\begin{proof}[Proof of Lemma~$\ref{lemma=LEFandRF}$]
Item $(1)$ follows from $(\star)$. The former-half of $(2)$ is proved in the viewpoint $(1)$ above because $\mathcal{Q}_{\mathbf{G}}$ exactly equals the set of all marked groups satisfying $r_1,\ldots, r_l$, where $r_1,\ldots, r_l$ determine the relations of $\mathbf{G}$. Then, the latter-half follows; see the remark above.
\end{proof}

In this paper, for a Cayley convergent sequence $(\mathbf{G}_m)_{m}$, $\mathbf{G}_m\stackrel{\mathrm{Cay}}{\to}\mathbf{G}$, we call it a \textit{LEF approximation} of $\mathbf{G}$ if $\mathbf{G}_m$ is finite for all $m$. We call it a \textit{LEA approximation} if $\mathbf{G}_m$ is amenable for all $m$. We call it an \textit{RF approximation} if $\mathbf{G}_m$ is of the form $(G/N_m;S\ \mathrm{mod}\ N_m)$ such that for every $m\in \mathbb{N}$, $N_m$ is a finite index normal subgroup of $G$ and $\mathrm{liminf}_{m\to \infty}N_m=\{e_G\}$ holds; see the remark above the proof of Lemma~\ref{lemma=LEFandRF}.
The latter-half of $(2)$ above says that a LEF approximation of a \textit{finitely presented} marked group is eventually an RF approximation. In the next subsection, we show examples of LEF approximations that are far from RF approximations; see also Subsections~\ref{subsection=SymmetricGroup} and \ref{subsection=ExtendedLinearGroups}.

The following lemma is quite easy to show; we omit the proof of it.

\begin{lemma}\label{lemma=LEFsubgroups}
Let $((L_m;v_1^{(m)},\ldots ,v_{\ell}^{(m)}))_{m\in \mathbb{N}}$ be a convergent sequence to $(L_{\infty};v_1^{(\infty)},\ldots ,v_{\ell}^{(\infty)})$ in the Cayley topology. Let $k\in \mathbb{N}_{\geq 1}$ and let $\omega_1,\ldots, \omega_k$ be words in $F_{\ell}$. For every $m\in \mathbb{N}\cup\{\infty\}$, let $s_j^{(m)}=\omega_j(v_1^{(m)},\ldots ,v_{\ell}^{(m)})$ for every $j\in [k]$. Let $G_m$ be the subgroup of $L_m$ generated by $s_j^{(m)}$, $j\in [k]$. Then $((G_m;s_1^{(m)},\ldots ,s_k^{(m)}))_{m\in \mathbb{N}}$ converges to $(G_{\infty};s_1^{(\infty)},\ldots ,s_k^{(\infty)})$ in the Cayley topology.
\end{lemma}

\subsection{Examples to apply Theorem~\ref{mthm=theoremA}}\label{subsection=Examples}
We here exhibit several examples that respectively express points $(a)$, $(b)$ and $(c)$ in Introduction.

\begin{example}[Point $(a)$]\label{example=Grigorchuk}
We describe an example of $(\mathbf{G}_m)_{m}$ whose Cayley boundary is a continuum consisting only of amenable groups; this example was brought to the authors by Grigorchuk. In \cite[Section~2]{Grigorchuk}, Grigorchuk defined $4$-marked groups $\mathbf{G}_{\omega}=(G_{\omega};a_{\omega},b_{\omega},c_{\omega},d_{\omega})$ associated with an (infinite) string $\omega$ in $\Omega=\{0,1,2\}^{\mathbb{N}}$. (For instance, for $\omega=(012)=012012\cdots$, $G_{\omega}$ is the first example of groups of intermediate growth.) Grigorchuk proved that for every $\omega \in \Omega$, $G_{\omega}$ is amenable and RF. Set $\Omega'$ to be the set of strings in $\omega\in \Omega$ that are not eventually constant (in other words, $\Omega'=\Omega\setminus \Omega_2$, where $\Omega_2$ is as in \cite{Grigorchuk}). We gather all finite marked group quotients  of $\mathbf{G}_{\omega}$ and take the union of such for all $\omega\in \Omega'$. Since every such marked group is finite, we can enumerate all the elements in this union (say, by order of the group). 
\begin{proposition}\label{proposition=Grigorchuk}
Let $(\mathbf{G}_m)_{m\in \mathbb{N}}$ be an enumeration of the above union of finite marked groups. Then $\partial_\mathrm{Cay}(\mathbf{G}_m)_{m}$ is a continuum of amenable marked groups. Moreover, it has a continuum of isomorphism classes of groups.
\end{proposition}
\begin{proof}
Grigorchuk  showed that $\partial_\mathrm{Cay}(\mathbf{G}_m)_{m}$ above coincides with the set $\{\tilde{\mathbf{G}}_{\omega}:\omega\in \Omega\}$; see \cite[Section~6 and Proposition~6.2]{Grigorchuk} for details. For $\omega\in \Omega'$, $\tilde{\mathbf{G}}_{\omega}=\mathbf{G}_{\omega}$; otherwise, $\tilde{G}_{\omega}$ is virtually solvable. The final assertion follows from \cite[Section~5]{Grigorchuk}.
\end{proof}

By $(i)$ of Theorem~\ref{mthm=theoremA}, $\coprod_{m\in \mathbb{N}}\mathrm{Cay}(\mathbf{G}_m)$ for $(\mathbf{G}_m)_{m\in \mathbb{N}}$ as in Proposition~\ref{proposition=Grigorchuk} has property A.
\end{example}

\begin{example}[Point $(b)$]\label{example=Point(b)}
We explain (a class of) groups that are LEF but not RF. Our first examples come from (restricted) \textit{$($permutational$)$ wreath products}. Let $G,H$ be  a group, $I$ be a set and $I\curvearrowleft H$ be a transitive action. Then the (restricted) \textit{permutational wreath products} $G\wr_I H$ is defined as the semidirect product $(\bigoplus_{I}G)\rtimes H$, where $H$ acts on $(\bigoplus_{I}G)$ by permutations of coordinates. If $I=H$ and the action $H\curvearrowleft H$ is by the right multiplications, then $G\wr_H H$ is usually written as $G\wr H$, and it is called the (restricted) \textit{wreath product} of $G$ and $H$. 

For $\mathbf{G}=(G;s_1,\ldots ,s_{k_1})\in \mathcal{G}(k_1)$ and $\mathbf{H}=(H;t_1,\ldots ,t_{k_2})\in \mathcal{G}(k_2)$, define a new marked group $\mathbf{G}\wr\mathbf{H} \in \mathcal{G}(k_1+k_2)$ as follows. The group is $G\wr H$, and the marking is $(s_1\delta_{e_H},e_H),(s_2\delta_{e_H},e_H),\ldots ,(s_{k_1}\delta_{e_H},e_H),(e, t_1),(e,t_2),\ldots ,(e,t_{k_1}))$. Here $e$ is the group unit of $\bigoplus_H G$, and for $g\in G$ and $h\in H$, $g\delta_{h} \in \bigoplus_H G$ denotes the element whose $h$-entry is $g$ and all of whose other entries are $e_G$.

\begin{proposition}\label{proposition=WreathProduct}
\begin{enumerate}[$(1)$]
  \item $(${\S}$2.4$. Theorem in \cite{VershikGordon}$)$ If $\mathbf{G}_m\stackrel{\mathrm{Cay}}{\to} \mathbf{G}_{\infty}\in \mathcal{G}(k_1)$ and if $\mathbf{H}_n\stackrel{\mathrm{Cay}}{\to} \mathbf{H}_{\infty}\in \mathcal{G}(k_2)$, then as $\min\{m,n\}\to \infty$,
\[
\mathbf{G}_m\wr\mathbf{H}_n\stackrel{\mathrm{Cay}}{\longrightarrow} \mathbf{G}_{\infty}\wr \mathbf{H}_{\infty}\quad \textrm{in}\quad \mathcal{G}(k_1+k_2).
\]
In particular, the LEF property and the LEA property are respectively closed under taking wreath products.
  \item $($Theorem~$3.2$ in \cite{Gruenberg}$)$ Let $H$ be an infinite  group. If $G\wr H$ is RF, then $G$ must be abelian.
\end{enumerate}
\end{proposition}
See \cite[Lemma~4.7 and Remark~4.8]{MimuraRF} for more details; these also provide another pedagogical example of the Cayley convergence other than Lemma~\ref{lemma=Tessera}.

Hence, $(\mathbb{Z}\wr \mathbb{Z})\wr \mathbb{Z}$ is one example of amenable LEF group that is not RF. Item $(1)$ in the proposition above should be clear to the readers who know a relationship between $G\wr H$ and random walks on $H$ with lamps in $G$. 

Another example is a (finitely generated) simple infinite amenable group. H. Matui \cite{Matui} and Juschenko--Monod \cite{JuschenkoMonod} provided the first example of such groups, the commutator subgroups of the topological full group of a minimal subshift on the Cantor set. In particular, these groups are not RF. On the other hand, Grigorchuk and K. Medynets \cite{GrigorchukMedynets} showed that they are LEF.

An example concerning amenable but not LEF groups can be obtained from work of L. Bartholdi and A. Erschler \cite{BartholdiErschler} on pre-orderings on the space of marked groups; see also \cite{Arzhantseva}. In \cite[Proposition~6.15]{BartholdiErschler}, they showed that if $G$ and $H$ are finitely generated group with $H$ infinite and if $G$ does not satisfy any identity (law), then there exist $k\in \mathbb{N}_{\geq 2}$ and a sequence of markings $(S_m)_{m\in\mathbb{N}}$ of $G\wr H$ such that $\lim_{m}(G\wr H;S_m)=\mathbf{F_k}$. We apply this for $G$ being $H_3$, one of \textit{Houghton groups}; see \cite[Subsection~5.3]{deCornulierGuyotPitsch}. The group $H_3$ is a \textit{finitely presented} group that admits a short exact sequence
\[
1\longrightarrow \mathrm{Sym}_{<\aleph_0}(\mathbb{Z}) \longrightarrow H_3 \longrightarrow \mathbb{Z}^2 \longrightarrow 1;
\]
recall our notation of symmetric groups from Introduction. Since both of $\mathrm{Sym}_{<\aleph_0}(\mathbb{Z})$ (this is a locally finite group) and $\mathbb{Z}^2$ are amenable, so is  $H_3$. Note that $\mathrm{Sym}_{<\aleph_0}(\mathbb{Z})$ contains a copy of $\mathrm{Alt}(\mathbb{Z})$ with index $2$. Here $\mathrm{Alt}(\mathbb{Z})$ denotes the alternating group over $\mathbb{Z}$, namely, the inductive limit of alternating groups on $[-m.m]$ over as $m\to \infty$. Since $\mathrm{Alt}(n)$ is simple for all $n\in \mathbb{N}_{\geq 5}$, the group $\mathrm{Alt}(\mathbb{Z})$ is simple. By Lemma~\ref{lemma=LEFandRF} and by simplicity of $\mathrm{Alt}(\mathbb{Z})$, $H_3$ is not LEF. Since $H_3$ contains a copy of $\mathrm{Sym}_{<\aleph_0}(\mathbb{Z})$, it does not satisfy any identity.

Therefore, by the aforementioned result in \cite{BartholdiErschler}, the group $H_3\wr\mathbb{Z}$, which is amenable, admits a sequence of markings $(S_m)_m$ with respect to which the marked groups converge to some (non-abelian) free marked group $\mathbf{F_k}$. By Theorem~\ref{mthm=theoremA}, the disjoint union $\bigsqcup_m\mathrm{Cay}(H_3\wr \mathbb{Z};S_m)$ does not have property A.

We here note the following: R. Willett \cite{Willett} showed that the (coarse) disjoint union of \textit{finite} connected graphs of every degree $\geq 3$ whose girth tends to $+\infty$ does not have property A. However, it can\textit{not} be generalized to the disjoint union of infinite graphs in the full generality: For instance, the disjoint union of (countable copies of) the infinite $4$-regular tree has property A because a single copy has property A. Therefore, our example concerning $H_3\wr\mathbb{Z}$ does not directly follow from Willett's result. In our case, these graphs are \textit{Cayley graphs} of \textit{amenable} groups, and we utilize it in the proof.

\end{example}

\begin{example}[Point $(c)$]\label{example=Point(c)}
We provide an example for which the phenomenon in $(c)$ occurs other than Example~\ref{example=SpecialLinear}. This is done by replacing the coefficient ring $\mathbb{F}_{p^{n_m}}$ of unimodular linear groups with some other rings; recall Remark~\ref{remark=OtherRings}. Fix $(l_m)_{m\in \mathbb{N}_{\geq 3}}$ a sequence of integers at least $2$ such that $\lim_{m\to \infty }l_m=+\infty$. Set $G_m=\mathrm{SL}^{\pm}(m,\mathbb{Z}/l_m\mathbb{Z})$ and take two markings $S_m$, $T_m$ as follows:
\begin{itemize}
  \item Set $S_m=(\sigma^{(m)},\delta^{(m)},\tau^{(m)})$, where $\sigma^{(m)}$ and $\tau^{(m)}$ are the matrices with exactly the same entries of $0$ and $1$ as in, respectively, $\sigma^{(m)}$ and $\tau^{(m)}$ in Example~\ref{example=SpecialLinear}. The element $\delta^{(m)}$ is the  matrix with exactly the same entries of $0$, $1$ and $-1$ as in $\delta^{(m)}$ in Example~\ref{example=SpecialLinear}.
  \item Set $T_m=(\sigma^{(m)},\sigma'^{(m)}, \delta^{(m)},\tau^{(m)})$, where $\sigma^{(m)}$, $\delta^{(m)}$ and $\tau^{(m)}$ are as in the item above, and $\sigma'^{(m)}={}^t\sigma^{(m)}$.
\end{itemize}

Similar to Corollary~\ref{mcor=ChoiceOfGenerators}, we will show the following in Subsection~\ref{subsection=ExtendedLinearGroups}.

\begin{proposition}\label{proposition=ChoiceOfGenerators}
Let $G_m$, $S_m$, $T_m$ be as in the current example for $m\in \mathbb{N}_{\geq 3}$.\begin{enumerate}[$(1)$]
    \item Then, $V=V_{(l_m)_m}=\bigsqcup_{m\geq 3} \mathrm{Cay}(G_m;S_m)$ has property A.
  \item Then, $W=W_{(l_m)_m}=\bigsqcup_{m\geq 3} \mathrm{Cay}(G_m;T_m)$ does not admit a coarse embedding into a Hilbert space.
\end{enumerate}
\end{proposition}
\end{example}
M. Kassabov pointed out to the authors that essentially the same sequence as $(\mathrm{Cay}(G_m;S_m))_{m\geq 3}$ of the current example was studied by A. Lubotzky and B. Weiss \cite{LubotzkyWeiss} in the context of expanders; see also examples in Remark~\ref{remark=SpecialLinearGroup}.

\subsection{Relation between the Cayley boundary and diagonal products of marked groups}\label{subsection=OtimesProduct}
In \cite[Definition~4.3]{KassabovPak}, Kassabov and I. Pak gave the definition of the diagonal product of marked groups as follows; in \cite{KassabovPak}, they used the terminology ``\textit{$\bigotimes$-products.}'' See also \cite[Subsection~2.1]{BrieusselZheng}.

\begin{definition}[Diagonal product]
Let $(\mathbf{G}_m=(G_m;s_1^{(m)},\ldots ,s_k^{(m)}))_{m\in\mathbb{N}}$ be a sequence in $\mathcal{G}(k)$. Then we define $\Delta_{m\in \mathbb{N}} (\mathbf{G}_m) \in \mathcal{G}(k)$ as follows: Set $\tilde{s}_1,\ldots ,\tilde{s}_k$ by $\tilde{s}_j=(s_j^{(0)},s_j^{(1)},s_j^{(2)},\ldots )\in \prod_{m\in \mathbb{N}} G_m$ for $j\in [k]$. Then, define $\Delta_{m\in \mathbb{N}} (\mathbf{G}_m) =(\Lambda;\tilde{s}_1,\ldots ,\tilde{s}_k)$, where $\Lambda$ is the subgroup of $\prod_{m\in \mathbb{N}} G_m$ generated by $\tilde{s}_1,\ldots ,\tilde{s}_k$.
\end{definition}

For every $n\in \mathbb{N}$, the projection onto $n$-th coordinate gives a marked group projection map $\Delta_{m\in \mathbb{N}} (\mathbf{G}_m)\twoheadrightarrow \mathbf{G}_n$. If $(\mathbf{G}_m)_m$ is an RF approximation of $\mathbf{G}$, then $\Delta_{m\in \mathbb{N}} (\mathbf{G}_m)$ coincides with $\mathbf{G}$. In what follows, we describe a relationship between this notion and the Cayley boundary of $(\mathbf{G}_m)_m$.

An \textit{ultrafilter} $\mathcal{U}$ over $\mathbb{N}$ has a one-to-one correspondence to a probability $\{0,1\}$-valued mean $\nu$ (finitely additive measure) that is defined on all subsets of $\mathbb{N}$: The corresponding $\mathcal{U}\longleftrightarrow \nu$ is 
\[
A\subseteq \mathbb{N} \textrm{ is in $\mathcal{U}$}\quad \Longleftrightarrow \quad \nu(A)=1.
\]
An ultrafilter over $\mathbb{N}$ is said to be \textit{non-principal} if it does not correspond to the Dirac mass at a point in $\mathbb{N}$; equivalenly, if it includes the \textit{cofinite filter} $\mathcal{U}_{\mathrm{cofin}}=\{A\subseteq \mathbb{N}: \sharp (\mathbb{N}\setminus A)<\infty\}$ (as a subfilter). 

For marked groups $(\mathbf{G}_m)_{m\in \mathbb{N}}$ and $\mathbf{H}$ in $\mathcal{G}(k)$ and for a non-principal ultrafilter $\mathcal{U}$, we write that $\lim_{\mathcal{U}}\mathbf{G}_m=\mathbf{H}$ if the following condition is satisfied: For every $R\in \mathbb{N}$, 
\[
\{m\in \mathbb{N}: B_{\mathrm{CayD}(\mathbf{G}_m)}(e_{G_m},R) \cong B_{\mathrm{CayD}(\mathbf{H})}(e_{H},R)\} \in \mathcal{U},
\]
where $\cong$ denotes the isomorphism as rooted diagrams. (Recall from $(\star)$ that the usual Cayley convergence amounts to requiring this with respect to $\mathcal{U}_{\mathrm{cofin}}$.) By compactness of $\mathcal{G}(k)$, for each choice of $\mathcal{U}$, $\lim_{\mathcal{U}}\mathbf{G}_m$ always exists, and the limit is unique for a fixed $\mathcal{U}$. If $(\mathbf{G}_m)_m$ is a convergence sequence, then $\lim_{\mathcal{U}}\mathbf{G}_m$ coincides with $\lim_{m\to \infty}{\mathbf{G}_m}$ for every non-principal ultrafilter $\mathcal{U}$. The Cayley boundary $\partial_{\mathrm{Cay}}((\mathbf{G}_m)_{m\in\mathbb{N}})$ exactly equals the set $\{\lim_{\mathcal{U}}\mathbf{G}_m$: $\mathcal{U}$ is a non-principal ultrafilter.$\}$.

Going back to the diagonal product, we set 
\[
N_{\mathcal{U}}=\left\{g=(g_0,g_1,\ldots )\in \prod_{m\in \mathbb{N}}G_m :\{m\in \mathbb{N}:g_m= e_{G_m}\}\in \mathcal{U}\right\}
\]
for each non-principal ultrafilter (for the cofinite filter, we can define $N_{\mathcal{U}_{\mathrm{cofin}}}$ in a similar manner, and it then equals $\bigoplus_{m\in \mathbb{N}}G_m$). Then $N_{\mathcal{U}}$ is a normal subgroup of $\prod_mG_m$, and we have a marked group isomorphism
\[
\lim_{\mathcal{U}}\mathbf{G}_m \cong (\Lambda/(N_{\mathcal{U}}\cap \Lambda); (\tilde{s}_1,\ldots, \tilde{s}_k)\ \mathrm{mod}\ (N_{\mathcal{U}}\cap \Lambda)).
\]
Compare with \cite[Lemma~4.6]{KassabovPak}.

In particular, for a LEA approximation $(\mathbf{G}_m)_m$ of an amenable group, the group $\Lambda$ as in $\Delta_m (\mathbf{G}_m)$ is amenable, as we also observed in our Part III paper \cite[Proposition~5.2]{MOSSPartIII}.

\begin{remark}\label{remark=FromLEFToRF}
If $(\mathbf{G}_m)_m$ is a LEF approximation of $\mathbf{G}_{\infty}$, then the underlying group $\Lambda$ of $\Delta_{m} (\mathbf{G}_m)$ is RF. To obtain a concrete RF approximation of $\Lambda$, for every $n\in \mathbb{N}$, consider 
\[
\mathbf{H}_n=\Delta_{m=0}^n (\mathbf{G}_m).
\]
Here the right-hand side is defined as $(H_n;\tilde{s}_1^{(n)},\ldots ,\tilde{s}_k^{(n)})$, where for every $j\in [k]$,
\[
\tilde{s}_j^{(n)}=(s_j^{(1)},\ldots ,s_j^{(n)}) \in \prod_{m=0}^n G_m.
\]
Then, $(\Lambda, s_1,\ldots ,s_k)$ admits a marked group quotient map onto each $\mathbf{H}_n$ via the  projection $\prod_{m\in \mathbb{N}}G_m \twoheadrightarrow \prod_{m=0}^n G_m$. Set $\tilde{N}_n$ be the kernel of the marked group quotient map. Then $(\tilde{N}_n)_{n\in \mathbb{N}}$ is a \textit{nested} sequence (recall the remark below Lemma~\ref{lemma=LEFandRF}) such that $\bigcap_{n\in \mathbb{N}} \tilde{N}_n=\{e_{\Lambda}\}$. Hence $(\mathbf{H}_n)_{n\in \mathbb{N}}$ is an \textit{RF} approximation of $(\Lambda;\tilde{s}_1,\ldots,\tilde{s}_k)$. 

By considering the sequence of marked groups $(\Delta_{m\in \mathbb{N}_{\geq n}}(\mathbf{G}_m))_{n\in \mathbb{N}}$, we observe that a LEF group is the Cayley limit of some RF groups. 
\end{remark}




\subsection{Large scale geometry: Property A, coarse embeddings and (coarse) disjoint union}\label{subsection=LargeScale}

In this subsection, we introduce basic terminologies in study of large scale geometry.
Let $X = (X, d)$ be a metric space.
In this paper, we often allow the value of the metric function $d$ to be $+\infty$. In that case, we call $X$ a \textit{generalized metric space} in the present paper.

\begin{definition}\label{DefinitionLocalFin}
The (generalized) metric space $(X, d)$ is said to be {\it uniformly locally finite} or to have {\it bounded geometry}, if
for each positive real number $R$, $\sup_{x \in X} \sharp ({B}_X(x,R)) < \infty$. Here ${B}_X(x,R)$ denotes the closed ball centered at $x$ of radius $R$.

For a sequence $(X_m, d_m)_m$ of (generalized) metric spaces, we say that $(X_m, d_m)_m$ is \textit{equi-uniformly locally finite} if the supremum above is uniformly finite on $m$ for every $R\in \mathbb{R}_{\geq0}$.
\end{definition}

Examples of our concern are Cayley graphs of $k$-marked groups and (coarse) disjoint unions of them (for fixed $k$). We here give the definition of the latter one(s).

\begin{definition}\label{definition=CoarseDisjointUnion}
Let $(X_m,d_m)_{m\in \mathbb{N}}$ be a sequence of metric spaces. 
\begin{enumerate}[$(1)$]
 \item The \textit{disjoint union} of $((X_m,d_m)_{m\in \mathbb{N}}$ is the generalized metric space constructed as follows: The space is $\bigsqcup_{m\in \mathbb{N}} X_m$, and we equip it with a generalized metric $d$ on that by $d(x,y)=d_m(x,y)$ if $x$ and $y$ are in the same component $X_m$; otherwise $d(x,y)=\infty$. We simply write the disjoint union (as a generalized metric space) as $\bigsqcup_{m\in \mathbb{N}} X_m$.
 \item Assume that each $X_m$, $m\in \mathbb{N}$ has a finite diameter (namely,  for every $m\in \mathbb{N}$, the maximal distance between two points in $X_m$ is finite). Then a \textit{coarse disjoint union} of $((X_m,d_m)_{m\in \mathbb{N}}$, which is a (genuine) metric space, is defined as follows: The underlying space is again $\bigsqcup_{m\in \mathbb{N}} X_m$. We endow it with a metric function $d'$ as follows.

\begin{eqnarray*}
d'(x,y)=
\left\{
\begin{array}{cll}
d_m (x, y), & \textrm{if}\quad (x, y) \in X_m \times X_m,\\
\mathrm{diam}(X_m) + \mathrm{diam}(X_n) + m + n, 
& \textrm{if}\quad (x, y) \in X_m \times X_n ,\ m \neq n,
\end{array}
\right.
\end{eqnarray*}
We write a coarse disjoint union as $\coprod_{m\in \mathbb{N}} X_m$.
\end{enumerate}
\end{definition}
For a sequence $(X_m)_m$ of equi-uniformly locally finite metric spaces, $\bigsqcup_m X_m$ is a uniformly locally finite generalized metric space. If $\mathrm{diam}(X_m)<\infty$ for all $m$, then $\coprod_m X_m$ is a uniformly locally finite metric space.

\begin{remark}\label{remark=CoarseDisjointUnion}
The specific form of this metric function $d'$ on the coarse disjoint union is not important; what is required is $\mathrm{dist}(X_m,X_n)\geq \max\{\mathrm{diam}(X_n),\mathrm{diam}(X_n)\}$ and $\mathrm{dist}(X_m,X_n)$ tends to $+\infty$ as $\min\{m,n\}$ tends to $\infty$ with satisfying $m\ne n$. Hence, unlike the disjoint union, there is ambiguity on concrete metrics of coarse disjoint unions, which is irrelevant in this paper. 

The notion of disjoint unions is natural in the framework of coarse spaces, but we will not go in details of that.
\end{remark}


We give one definition of property A of Yu.

\begin{definition}[Definition 11.35 of \cite{RoeLectureNote}]
\label{Lemma Property A Hulanicki}
A uniformly locally finite (generalized) metric space $(X, d)$ has \textit{property A} if the following condition holds:
For every $C > 0$ and for every (small) $\epsilon > 0$, 
there exists a map $X \ni x \mapsto \eta_x \in \ell_2 (X)$ assigning unit vectors such that
\begin{itemize}
\item
$\sup \{d (x, y) : x \in \mathrm{supp}(\eta_y)\} < \infty$,
\item
if $d(x, y) < C$, then $\|\eta_x - \eta_y \| < \epsilon$.
\end{itemize}
\end{definition}

\begin{remark}\label{Remark Hilbert}
We may replace $\ell_2 (X)$ with any other infinite dimensional Hilbert space $\mathcal{H}$ 
(following an idea in \cite[Theorem 3]{PaperBrodzkiNibloWright}).
When we use $\mathcal{H}$, the first condition of $\eta$ is replaced with $\sup\{ d(x, y) : \langle \eta_x, \eta_y \rangle \neq 0\}<\infty$.
As a corollary, property A passes to subsets.
\end{remark}

We exhibit some examples of metric spaces with property A. For finitely generated groups (viewed as metric spaces by some word length/Cayley graph), amenable groups and linear groups \cite{GuentnerHigsonWeinberger} both have property A. N. Ozawa \cite[Corollary 3]{BoundaryAme} showed that groups hyperbolic relative to (finitely many) subgroups with property A has property A. Quite recently, it is showed that $\mathrm{Out}(F_n)$ have property A \cite{BestvinaGuirardelHorbez} for all $n$. A  uniformly locally finite metric space with finite \textit{asymptotic dimension} \cite[Section 1E]{GromovAsymptotic} has property A \cite[Lemma 4.3]{HigsonRoe}. On the other hand, there exist some groups without property A; see \cite{Gromovrandomwalk}, \cite{ArzhantsevaDelzant}, \cite{ArzhantsevaOsajda} and \cite{Osajda}. D. Osajda \cite{OsajdaRF}, building on earlier work of \cite{Osajda} and \cite{ArzhantsevaOsajda}, constructed an RF group without property A. We employ these constructions in our Part II paper; see \cite[Theorem~C and Subsection~9.2]{MSPartII}.

%




We recall a  notion that amounts to an inclusion in large scale geometry, which is called a coarse embedding.

\begin{definition}\label{definition=CoarseEmbedding}
Let $(X,d_X)$ and $(M,d_M)$ be  generalized metric spaces. 
\begin{enumerate}[$(1)$]
 \item A (possibly discontinuous and possibly non-injective) map $f \colon X \rightarrow M$ is called a \textit{coarse embedding} if there exist non-decreasing functions $\rho, \omega \colon [0, \infty) \rightarrow [0, \infty)$ that are proper (namely, $\lim_{r\to +\infty}\rho(r)=\lim_{r\to +\infty}\omega(r)=+\infty$) such that for all $(x_1,x_2)\in X\times X$ with $d_X(x_1,x_2)<\infty$, it holds that
\[
\rho(d_X(x_1, x_2)) \le d_M(f(x_1), f(x_2)) \le \omega(d_X(x_1, x_2)).
\]
(In particular, if $d_X(x_1,x_2)<\infty$, then $d_M(f(x_1),f(x_2))<\infty$.)
 \item In the definition above, we call $\rho$, $\omega$ and $(\rho,\omega)$, respectively, a \textit{compression function}, an \textit{expansion function} and a \textit{control pair} for $f$. A pair $(\rho,\omega)$ of non-decreasing proper functions $[0, \infty) \rightarrow [0, \infty)$ is called a \textit{control pair} for $X$ into $M$ if there exists a coarse embedding $f\colon X\to M$ for which $(\rho,\omega)$ is a control pair. Define $\mathcal{CP}_M(X)$ to be the the set of all control pairs for $X$ into $M$.
\end{enumerate}
\end{definition}

For a (generalized) metric space $X$ that has property A, there exists a coarse embedding of $X$ into a Hilbert space; see \cite{Yu}.

For a sequence $(X_m=(X_m,d_m))_m$ of equi-uniformly locally finite metric spaces each of whose diameter is finite, it is straightforward to see that the disjoint union $\bigsqcup_m X_m$ admits a coarse embedding into a Hilbert space if and only if so does a coarse disjoint union $\coprod_mX_m$. Indeed, to construct a coarse embedding into a Hilbert space from $\coprod_mX_m$ from one from $\bigsqcup_m X_m$, we send each $X_m$ to an appropriate codimension $1$ subspace $\mathcal{H}_m$ of a Hilbert space in such a way that $\mathrm{dist}(\mathcal{H}_m,\mathcal{H}_n)$ is sufficiently big in terms of $(m,n)$ with $m\ne n$. Concerning on property A, we prove the following.

\begin{proposition}\label{PropositionA}
Let $(X_m)_{m\in \mathbb{N}}$ be a sequence of equi-uniformly locally finite metric spaces such that $\mathrm{diam}(X_m)<\infty$ for all $m\in \mathbb{N}$. Then, the disjoint union $\bigsqcup_{m \in \mathbb{N}} X_m$ has property A if and only if so does a coarse disjoint union $\coprod_{m \in \mathbb{N}} X_m$.
\end{proposition}

This proposition almost immediately follows from a general statement on \textit{union permanence} of certain coarse geometric properties; see \cite[Proposition~2.20]{Guentner}. However, for the reader's convenience, we include a direct proof of Proposition~\ref{PropositionA}, as follows.

\begin{proof}
We prove this for the specific choice of the metric $d'$ on $\coprod_m X_m$ as in Definition~\ref{definition=CoarseDisjointUnion}; the proof works for all the other choices that satisfy the conditions in Remark~\ref{remark=CoarseDisjointUnion}.  Let $d$ be the metric function on the disjoint union $X = \bigsqcup_{m\in \mathbb{N}} X_m$.

Assume that the disjoint union $X$ has property A.
For every $\epsilon > 0$, and $C > 0$, there exists a map $X \ni x \mapsto \eta_x \in \ell_2 (X)$ assigning unit vectors such that 
\begin{itemize}
\item
if $d(x, y) < C$, then $\|\eta_x - \eta_y\| < \epsilon$,
\item
there exists $R > 0$ such that $C < R$ and 
$d(x, y) < R$ whenever $x \in \mathrm{supp}(\eta_y)$.
\end{itemize}
Choose and fix some unit vector $\eta_0$ in $\ell_2 \left( X_1 \sqcup \cdots \sqcup X_{\lfloor R\rfloor} \right)$, where $\lfloor R\rfloor$ denotes the integer part of $R$.
Define $\tilde\eta \colon \coprod_{m\in \mathbb{N}} X_m \to \ell_2(\coprod_{m\in \mathbb{N}} X_m))$ as
\begin{eqnarray*}
\left\{
\begin{array}{cccc}
\tilde\eta_x = \eta_0, &\quad \textrm{if}\quad x \in X_m \textrm{ with } m \le R, \\
\tilde\eta_x = \eta_x, &\quad \textrm{if}\quad x \in X_m \textrm{ with }  m > R.\\
\end{array}
\right.
\end{eqnarray*}
If $x, y \in \coprod_{m\in \mathbb{N}} X_m$ and $x \in \mathrm{supp} (\tilde\eta_y)$, then
$(x, y) \in \left( X_1 \sqcup \cdots \sqcup X_{\lfloor R\rfloor} \right)^2$ or
$d'(x, y) \le d(x, y) \le R$.
It follows that 
\[\sup \{d'(x, y) : x \in \mathrm{supp}(\tilde\eta_y)\} < \infty.\]
Since $C < R$,
for $x, y \in \coprod_{m\in \mathbb{N}} X$,
if $d'(x, y) < C$, then $x, y \in X_1 \sqcup \cdots \sqcup X_{\lfloor R\rfloor}$ or $d(x, y) < C$.
It follows that if $d'(x, y) < C$, then $\| \tilde\eta_x - \tilde\eta_y \| < \epsilon$.
We hence conclude that the coarse disjoint union $\coprod_{m \in \mathbb{N}} X_m$ has property A.

The converse is straightforward.
\end{proof}

\section{Property A and amenability}
\label{section=Amenability}
In this section, we prove item $(i)$ of Theorem~\ref{mthm=theoremA}; in fact, we prove the following.

\begin{theorem}\label{theorem=Item(i)}
Let $(\mathbf{G}_m)_{m\in \mathbb{N}}$ be a sequence consisting of $\mathrm{amenable}$ marked groups in $\mathcal{G}(k)$. Then, the following are equivalent:
\begin{enumerate}[$(1)$]
  \item The disjoint union $\bigsqcup_m\mathrm{Cay}(\mathbf{G}_m)$ has property A.
  \item $($Under the condition that all $\mathbf{G}_m$, $m\in \mathbb{N}$, are finite$)$ A coarse disjoint union $\coprod_m\mathrm{Cay}(\mathbf{G}_m)$ has property A.
  \item The Cayley boundary $\partial_{\mathrm{Cay}}((\mathbf{G}_m)_m)$ is $\mathrm{pointwise}$ $\mathrm{amenable}$, that means, every element in $\partial_{\mathrm{Cay}}((\mathbf{G}_m)_m)$ is amenable.
  \item The Cayley boundary $\partial_{\mathrm{Cay}}((\mathbf{G}_m)_m)$ is $\mathrm{uniformly}$ $\mathrm{amenable}$ in the sense of Definition~$\ref{definition=UniformlyAmenable}$.
  \item The Cayley closure $\overline{(\mathbf{G}_m)_m}^{\mathrm{Cay}}$ is pointwise amenable.
  \item The Cayley closure $\overline{(\mathbf{G}_m)_m}^{\mathrm{Cay}}$ is uniformly amenable.
\end{enumerate}
\end{theorem}

\subsection{Uniform amenability}
Amenability of groups has considerably many equivalent characterizations. One of them for finitely generated groups is the \textit{F{\o}lner property}. 

\begin{definition}\label{defintion=Folner}
For $R\in \mathbb{N}$,
define a function $\mathrm{Rel}(\ \cdot\ , R)$ on
marked groups by
\begin{eqnarray*}
\mathrm{Rel}(\mathbf{G} , R)
=
\min
\left\{
\frac{\sharp(\partial_\mathbf{G}(Y,1))}{\sharp(Y)}\ 
:\ 
\emptyset \neq Y \subseteq {B}_{\mathbf{G}} (e_G, R)
\right\}.
\end{eqnarray*}
\end{definition}

By definition, $\mathrm{Rel}(\cdot, R)$ is non-increasing on $R\in \mathbb{N}$. The \textit{F{\o}lner set characterization} of amenability states that a finitely generated group $G$ is \textit{amenable} if and only if $\inf_R \mathrm{Rel}(\mathbf{G}, R) = 0$ for some (equivalently, all) marking $\mathbf{G}$ of $G$. 

\begin{definition}\label{definition=UniformlyAmenable}
A non-empty subset $K\subseteq \mathcal{G}(k)$ is said to be \textit{uniformly amenable} if 
\[
\inf_{R\in \mathbb{N}} \sup_{\mathbf{G}\in K}\mathrm{Rel}(\mathbf{G}, R) = 0
\]
holds true.
\end{definition}

\subsection{Topological properties of the set of amenable groups}
The goal of this subsection is to show that \textit{the uniformity is automatic} for a \textit{compact} subset of $\mathcal{G}(k)$ concerning amenability. More precisely, we will prove the following.

\begin{proposition}\label{propositon=AutomaticUniformity}
Let $K$ be a non-empty $\mathrm{compact}$ subset of $\mathcal{G}(k)$. Then the following are equivalent:
\begin{enumerate}[$(1)$]
  \item The set $K$ is pointwise amenable.
  \item The set $K$ is uniformly amenable.
\end{enumerate}
\end{proposition}

\begin{proof}
From $(\star)$ (and Lemma~\ref{lemma=Neighborhood}), the following can be observed (employ $\beta_{\mathbf{H},\mathbf{G},R+1}$):

\begin{lemma}\label{lemma=Continuous}
For every $R\in \mathbb{N}$, $\mathrm{Rel}(\ \cdot\ , R) \colon \mathcal{G}(k) \to \mathbb{R}_{\geq 0}$ is locally constant. In particular, it is continuous.
\end{lemma}

We only need to show [$(1)\Rightarrow (2)$]. The non-increasing sequence $(\mathrm{Rel}(\ \cdot\ , R))_{R\in \mathbb{N}}$ pointwise converges to $0$ on $K$ by the assumption $(1)$. By Lemma~\ref{lemma=Continuous}, the Dini Theorem shows that this convergence is uniform on $K$.
\end{proof}




We remark the following set-theoretic observation on amenable marked groups.

\begin{proposition}
\label{PropositionGdelta}
For every $k\in \mathbb{N}_{\geq 1}$, the set of amenable marked groups $\mathcal{A}(k) \subseteq \mathcal{G}(k)$ is an intersection of countably many open subsets of $\mathcal{G}(k)$. The relative topology on $\mathcal{A}(k)$ with respect to
the Cayley topology is metrizable by a complete metric.
\end{proposition}

\begin{proof}
For each $n\in \mathbb{N}$, define an open subset $\mathcal{F}_n$ of $\mathcal{G}(k)$ 
by 
\[\mathcal{F}_n = \bigcup_{R\in \mathbb{N}_{\geq 1}} 
\{\mathbf{G} \in \mathcal{G}(k) : \mathrm{Rel}(\mathbf{G}, R) < 2^{-n}\}.
\]
Then, it follows that
$\mathcal{A}(k)$ is identical to the intersection $\bigcap_{n\in \mathbb{N}_{\geq 1}} \mathcal{F}_n$.
The second assertion follows from a general theorem for Polish spaces (see, for instance, \cite[Theorem 3.11]{Kechris}).
\end{proof}

\subsection{Proof of Theorem~\ref{theorem=Item(i)}}
In this subsection, we will prove Theorem~\ref{theorem=Item(i)}. Our proof employs the concept of \textit{positive definite functions} on groups and \textit{GNS constructions} of them in an essential way. The readers who are not familiar with these topics may consult \cite[Appendix~C]{BookBekkadelaHarpeValette}.

We will utilize the following  three lemmata.

\begin{lemma}
\label{LemmaContraction}
Let $\varphi\colon \mathbf{F_k}\twoheadrightarrow \mathbf{G}$ be a marked group quotient map. Let $N$ be $\mathrm{Ker}\varphi \trianglelefteq F_k$.
Then for every $f_1, f_2 \in F_k$, 
\[
d_{\mathbf{G}}(\varphi(f_1),\varphi(f_2))=
\mathrm{dist}_{\mathbf{F_k}}(f_1f_2^{-1},N)
(=\min\{ d_{\mathbf{F_k}}(f_1 f_2^{-1}, x) : x \in N \}).
\]
\end{lemma}

\begin{proof}
Straightforward.
\end{proof}

\begin{lemma}
\label{LemmaPullBack}
Let $\varphi \colon F \rightarrow G$ be a homomorphism between groups.
Let $\phi$ be a positive definite function on $G$.
Then the function $\psi = \phi \circ \varphi$ on $F$ is positive definite.
\end{lemma}

\begin{proof}
By definition.
\end{proof}

\begin{lemma}
\label{LemmaPushForward}
Let $\varphi \colon F \rightarrow G$ be a surjective 
homomorphism between groups.
Denote by $N$ the kernel of $\varphi$.
Let $\psi$ be a positive definite function on $F$.
Suppose that $\psi(h) = 1$ for every $h \in N$.
Then 
there exists a positive definite function $\phi$ on $G$ such that
$\psi = \phi \circ \varphi$.
\end{lemma}

\begin{proof}
Let $(v, \mathcal{H}, \eta)$ be the GNS-construction
for the positive definite function $\psi$.
Then the function $\psi$ is expressed as
$\psi(f) = \langle v(f)\eta, \eta \rangle$.
Since $e_{F} \in N$,
we have that $\|\eta\|^2 = \psi(1_{F_k}) = 1$.
For every $h \in N$ and $f \in F$,
since $f^{-1} h f \in N$, it follows that 
\begin{eqnarray*}
\|v(h) v(f)\eta - v(f) \eta\|^2 
&=& 2 \| \eta \|^2 - 
2 \mathrm{Re}(\langle v(h f)\eta, v(f) \eta \rangle)\\
&=& 2 - 2\mathrm{Re}(\psi(f ^{-1} h f)) \\
&=& 0.
\end{eqnarray*}
Hence $N$ is included in the kernel of the group homomorphism $v$.
Hence, 
there exists
a group homomorphism $u \colon G \rightarrow \mathcal{U}(\mathcal{H})$ such that $v = u \circ\varphi$.
The function $\phi(g) = \langle u(g)\eta, \eta \rangle$
on $G$ satisfies the required conditions.
\end{proof}

\begin{proof}[Proof of Theorem~$\ref{theorem=Item(i)}$]
The equivalence [$(1)\Leftrightarrow (2)$] is by Proposition~\ref{PropositionA}. The equivalences [$(3)\Leftrightarrow (4)$] and [$(5)\Leftrightarrow (6)$] both follow from Proposition~\ref{propositon=AutomaticUniformity}. The equivalence [$(3)\Leftrightarrow (5)$] is by assumption. Hence, it suffices to show [$(1)\Rightarrow (3)$] and [$(6)\Rightarrow (1)$]. 

First, we demonstrate that $(6)$ implies $(1)$. 
Take arbitrary (big) $C > 0$ and (small) $\epsilon > 0$. Then by the assumption of $(6)$, 
there exists $R \in \mathbb{N}_{\geq 1}$ such that 
\[
\mathrm{Rel}(\mathbf{G}_m, R) < \frac{\epsilon^2}{C}, \quad \textrm{for all }m \in \mathbb{N}.
\]
It follows that for each $m$, there exists a finite subset $Y_m$ 
of $G_m$ satisfying
\begin{eqnarray*}
\sharp \left( \partial_{\mathbf{G}_m} \left( Y_m ,1\right) \right)
<
\frac{\epsilon^2}{C} \sharp \left( Y_m \right) \quad \textrm{and} \quad
Y_m \subseteq {B}_{\mathbf{G}_m} \left( e_{G_m},R \right).
\end{eqnarray*}
Define a unit vector $\eta_m$ in $\ell_2(G_m) \subseteq \ell_2(\bigsqcup_{m\in \mathbb{N}} G_m)$ by
\[
\eta_m = \frac{1}{\sqrt{\sharp(Y_m)}} 
\sum_{y \in (Y_m)^{-1}} \delta_y.
\]
For $g \in G_m$, we define a unit vector $\eta_{m,g} \in \ell_2(G_m)$
by
\[
\eta_{m,g} = \eta_m * \delta_g = \frac{1}{\sqrt{\sharp(Y_m)}} 
\sum_{y \in (Y_m)^{-1}} \delta_{y g}.
\]
The support of $\eta_{m,g}$ is included in 
${B}_{\mathbf{G}_m} (g,R) \subseteq G_m$.

For every pair $(g, h) \in (G_m)^2$, we have that
\begin{eqnarray*}
\| \eta_{m,g} - \eta_{m,h} \|^2
&=& \frac{\sharp((Y_m)^{-1} g \bigtriangleup (Y_m)^{-1} h)}
{\sharp(Y_m)}
= \frac{\sharp( g h^{-1} Y_m \bigtriangleup Y_m)}
{\sharp(Y_m)}\\
&\le&
d_{\mathbf{G}_m}(e_{G_m},gh^{-1}) \cdot
\frac{\sharp(\partial_{\mathbf{G}_m}(Y_m,1))}{\sharp(Y_m)}.
\end{eqnarray*}
For every $m \in \mathbb{N}$ and $g,h \in G_m$ with $d_{\mathbf{G}_m}(g, h) \le C$, we have that
\[
\| \eta_{m,g} - \eta_{m,h} \|^2
\le C
\frac{\sharp(\partial_{\mathbf{G}_m}(Y_m,1))}{\sharp(Y_m)}
< \epsilon^2.\]
It follows that the metric space $\bigsqcup_{m\in \mathbb{N}} \mathrm{Cay}(\mathbf{G}_m)$ has property A.

Finally, we prove the implication [$(1)\Rightarrow (3)$]. In what follows, we employ the limiting and averaging technique described in 
\cite[Proposition 11.39]{RoeLectureNote}.
Suppose that
$\bigsqcup_{m \in \mathbb{N}} \mathrm{Cay}(\mathbf{G}_m)$ has property A.
Take an  element $\mathbf{G}_{\infty}\in \partial_{\mathrm{Cay}}(\mathbf{G}_m)_m$ arbitrarily, and fix it. 
Take a subsequence of $(\mathbf{G}_m)_{m\in \mathbb{N}}$
that converges to $G$. 
Since property A passes to subsets, we may assume that
$(\mathbf{G}_m)_{m\in \mathbb{N}}$ converges to $\mathbf{G}_{\infty}$.
Let $\epsilon$ be an arbitrary (small) strictly positive number and let $C$ be an arbitrary (large) natural number.
Then there exist a natural number $R$ and a map 
\[
\eta :
\bigsqcup_{m\in \mathbb{N}} G_m 
\rightarrow \ell_2 \left( \bigsqcup_{m \in \mathbb{N}} G_m \right)
\]
assigning unit vectors and satisfying the following properties:
\begin{itemize}
\item
If $g, h \in G_m$ and $d_{\mathbf{G}_m}(g, h) \le C$, 
then $\|\eta_g - \eta_h\| < \epsilon$. 
\item
If $g \in G_m$, then $\mathrm{supp}(\eta_g) \subset 
{B}_{\mathbf{G}_m}(g,R) \subseteq G_m$.
\end{itemize}
We define a positive definite kernel $\Phi$ on $\bigsqcup_{m \in \mathbb{N}} G_m$ by
\[
\Phi(g, h) = \langle \eta_g, \eta_h \rangle, \quad
(g, h) \in \bigsqcup_{m \in \mathbb{N}} G_m \times \bigsqcup_{m \in \mathbb{N}} G_m.
\]
The kernel $\Phi$ satisfies the following:
\begin{itemize}
\item
For $g \in G_m$, $\Phi(g, g) = 1$.
\item
For $(g, h) \in G_m \times G_m$, 
if $d_{\mathbf{G}_m}(g, h) \le C$, then $|1 - \Phi(g, h)| < \epsilon$.
\item
For $(g, h) \in G_m \times G_m$, 
if $d_{\mathbf{G}_m}(g, h) > 2R$, then $\Phi(g, h) = 0$.
\end{itemize}

Since the group $G_m$ is amenable, there exists a left invariant mean $\mu_m$ on $G_m$.
For $g \in G_m$, define a bounded function $\Phi_g$ on $G_m$ by $\Phi_g(h) = \Phi(g^{-1} h, h)$. 
Define a function $\phi_m$ on $G_m$ by
\[ \phi_m(g) = \mu_m (\Phi_g). \]
The function $\phi_m$ satisfies the following:
\begin{itemize}
\item
If $d_{\mathbf{G}_m}(e_{G_m},g) \leq C$,
then $|1 - \phi_m(g)| < \epsilon$.
\item
If $d_{\mathbf{G}_m}(e_{G_m},g)>2R$,
then $\phi_m(g) = 0$.
\end{itemize}
We claim that $\phi_m$ is a positive definite function.
Let $\lambda$ denote the left translation action on $\ell_\infty(G_m)$:
$[\lambda_g (\zeta)](h) = \zeta(g^{-1} h)$.
For $n \in \mathbb{N}_{\geq 1}$, $g_1, g_2, \cdots, g_n \in G_m$
and $\alpha_1, \alpha_2, \cdots, \alpha_n \in \mathbb{C}$,
$\sum_{i, j = 1}^n 
\overline{\alpha_j} \alpha_i
\lambda_{g_j} \left( \Phi_{g_j^{-1} g_i} \right)$
is non-negative because
\[
\sum_{i, j = 1}^n 
\overline{\alpha_j} \alpha_i
\lambda_{g_j} \left( \Phi_{g_j^{-1} g_i} \right) (h)
=
\sum_{i, j = 1}^n 
\overline{\alpha_j} \alpha_i
\Phi(g_i^{-1} h, g_j^{-1} h).
\]
It follows that
\begin{eqnarray*}
\sum_{i, j = 1}^n 
\overline{\alpha_j} \alpha_i \phi_m(g_j^{-1} g_i)
=
\mu^{(m)}
\left( 
\sum_{i, j = 1}^n 
\overline{\alpha_j} \alpha_i
\lambda_{g_j} \left( \Phi_{g_j^{-1} g_i} \right)
\right)
\ge 0.
\end{eqnarray*}
Hence we prove the claim above.

Let $\psi_m$ be the composition
of the marked  quotient map $\varphi_m \colon \mathbf{F_k} \rightarrow \mathbf{G}_m$
and the function $\phi_m$. 
By Lemma \ref{LemmaPullBack}, 
$\psi_m$ is a positive definite function of $F_k$.
The function $\psi_m$ satisfies the following:
\begin{enumerate}
\item\label{condition1}
For $f \in F_k$, if $d_{\mathbf{F_k}}(e_{F_k},f) \le C$,
then $|1 - \psi_m(f)| < \epsilon$.
\item\label{condition2}
If $\mathrm{dist}_{\mathbf{F_k}}(f, \mathrm{Ker}(\varphi_m)) > 2R$, then $\psi_m(f) = 0$.
\end{enumerate}
Indeed, the second assertion follows from Lemma~\ref{LemmaContraction}.
Since $\psi_m(e_{F_k}) = 1$, we have 
$|\psi_m (f) | \le 1$ for every $f \in F_k$.
By the compactness of the disk 
$\{z \in \mathbb{C} : |z| \le 1\}$,
there exists a subsequence $(m(l))_{l\in \mathbb{N}}$
such that for every $f \in F_k$,
the sequence $(\psi_{m(l)}(f))_{l\in \mathbb{N}}$ converges.
Its limit $\psi_\infty(f)$ is also a positive definite function on $F_k$.
If $f$ is in the kernel of the marked quotient map $\varphi_{\infty} \colon \mathbf{F_k} \rightarrow \mathbf{G}_{\infty}$, then there exists $m_f\in \mathbb{N}$ such that $f$ is in the kernel of $\varphi_m \colon \mathbf{F_k} \rightarrow \mathbf{G}_m$ for all $m\geq m_f$; recall the original definition of the Cayley topology from Definition~\ref{definition=CayleyTopology}. 
It hence follows that $\psi_\infty$ is $1$ on the normal subgroup $\mathrm{Ker}(\varphi_{`\infty}) \trianglelefteq F_k$.
Lemma \ref{LemmaPushForward} shows that $\psi_\infty$
induces a positive definite function $\phi$ on $G_{\infty}$.

For $g \in G_{\infty}$, if $d_{\mathbf{G}_{\infty}}(e_{G_{\infty}},g)\leq C$, 
then $|1 - \phi(g)| < \epsilon$ by condition (\ref{condition1}).
By condition (\ref{condition2}),
for $f \in F_k$, if $\mathrm{dist}_{\mathbf{F_k}}(f, \mathrm{Ker}(\varphi_{\infty})) > 2R$, then $\psi_\infty(f) = 0$.
By Lemma~\ref{LemmaContraction}, it follows that if $d_{\mathbf{G}_{\infty}}(e_{G_{\infty}},g)> 2R$,
then $\phi(g) = 0$.

By varying $(C,\epsilon)$ appropriately (with according changes of $R=R(C,\epsilon)$), we may thus construct a sequence of positive definite functions $(\phi_n)_{n\in \mathbb{N}}$ on $G_{\infty}$ with $\phi_n(e_{G_{\infty}})=1$ such that all of $\phi_n$ are finitely supported and that $\phi_n$ converges to the constant $1$ function pointwise. It characterizes amenability of $G_{\infty}$. 
\end{proof}

\begin{proof}[Proof of $(i)$ of Theorem~$\ref{mthm=theoremA}$]
Immediate from [$(1)\Leftrightarrow (2) \Leftrightarrow (3)$] in Theorem~\ref{theorem=Item(i)}.
\end{proof}

\section{A-T-menability and coarse embeddings into a Hilbert space}\label{section=a-T-menable}
In this section, we will prove item $(ii)$ of Theorem~\ref{mthm=theoremA}. A finitely generated  group $G$ is said to be \textit{a-$\mathrm{T}$-menable} if there exists a coarse embedding $f$ from $\mathbf{G}$ to a Hilbert space $\mathcal{H}$ that is of the form 
\[
f(g)=\xi \cdot \alpha(g) , \quad \textrm{for all $g\in G$}
\]
for some $\xi\in \mathcal{H}$ and for some affine isometric action $\alpha\colon \mathcal{H} \curvearrowleft G$. Here we consdier a \textit{right} action because we equip  marked groups with \textit{right-invariant} metrics. (This definition does not depend on the choice of markings of $G$; however, the control pair for such an embedding will depend on it by multiplicative constants.) More shortly, $G$ is a-$\mathrm{T}$-menable if there exists a coarse embedding into a Hilbert space that is ``\textit{$G$-equivariant}.'' Through Schoenberg's argument, this condition is equivalent to the following: There exists a sequence of positive definite functions $(\phi_m \colon G \to \mathbb{R})_{m\in \mathbb{N}}$ which satisfies $\phi_m(e_G)=1$ and which is $c_0$ (namely, vanishing at infinity) such that $\phi_m$ pointwise converges to the constant $1$ function. (This property is called the \textit{Haagerup property}.) For instance, (finitely generated) amenable groups and  $F_2$ are a-$\mathrm{T}$-menable. For comprehensive treatises on a-$\mathrm{T}$-menable groups, see \cite{CCJJV} and \cite{BookBekkadelaHarpeValette}.

\begin{theorem}\label{theorem=a-T-menable}
Let $(\mathbf{G}_m)_{m\in \mathbb{N}}$ be a convergent sequence consisting of finite marked groups in $\mathcal{G}(k)$. Let $\mathbf{G}_{\infty}=(G_{\infty};S_{\infty})$ be the limit.
If $X=\bigsqcup_{m\in \mathbb{N} } \mathrm{Cay}(\mathbf{G}_m)$ admits a coarse embedding into a Hilbert space, then $G_{\infty}$ is a-$\mathrm{T}$-menable.
\end{theorem}

\begin{proof}
By  \cite[Theorem 5.2.8]{BookNowakYu},
for every (small) $\epsilon > 0$ and every (large) $C \in \mathbb{N}$, there exist
a Hilbert space $\mathcal{H}$,
a map $\eta_{\cdot} : X \rightarrow \mathcal{H}$,
and a non-decreasing function $\tilde{\rho} : \mathbb{N} \rightarrow [0, \sqrt{2})$ converging to $\sqrt{2}$ as $R\in \mathbb{N}$ tends to $+\infty$
satisfying that
\begin{itemize}
\item
for every $m$ and for every $x \in G_m$,
$\|\eta_x\| = 1$,
\item
for every $(x, y) \in (G_m)^2$,
$\| \eta_x - \eta_y \| \ge \tilde{\rho}(d_{\mathbf{G}_m}(x, y))$,
\item
if $d_{\mathbf{G}_m}(x, y) \le C$, then $\| \eta_x - \eta_y \| < \epsilon$.
\end{itemize}
Define a positive definite kernel $\Phi_m$ on $(G_m)^2$ by
$
\Phi_m(x, y) = \mathrm{Re} \langle \eta_x, \eta_y \rangle.
$
Then $\Phi_m$ satisfies
\begin{itemize}
\item
$\Phi_m(x, y) = (2 - \| \eta_x - \eta_y \|^2)/2 
\le (2 - \tilde{\rho}(d_{\mathbf{G}_m}(x, y))^2)/2$,
\item
if $d_{\mathbf{G}_m}(x, y) \le C$, then $\Phi_m(x, y) 
= 1 - \mathrm{Re} \langle \eta_x, \eta_x - \eta_y \rangle
\ge 1 - \| \eta_x - \eta_y \|
> 1 - \epsilon
$.
\end{itemize}
Define a function $\phi_m$ on $G_m$ by 
\[
\phi_m(g) = \frac{1}{\sharp(G_m)}
\sum_{h \in G_m} \Phi_m(g^{-1} h, h).
\]
For every $n \in \mathbb{N}$, $g_1$, $\cdots$, $g_n \in G_m$, and $\alpha_1$, $\cdots$, $\alpha_n \in \mathbb{C}$,
we have
\begin{eqnarray*}
\sum_{i, j = 1}^n \overline{\alpha_j} \alpha_i \phi_m(g_j^{-1} g_i)
&=&
\frac{1}{\sharp(G_m)}
\sum_{i, j = 1}^n \sum_{h \in G_m} 
\overline{\alpha_j} \alpha_i
\Phi_m(g_i^{-1} g_j h, h)\\
&=&
\frac{1}{\sharp(G_m)}
\sum_{h \in G_m}
\sum_{i, j = 1}^n 
\overline{\alpha_j} \alpha_i
\Phi_m(g_i^{-1} h, g_j^{-1} h)\\
&\ge& 0.
\end{eqnarray*}
It follows that $\phi_m$ is a positive definite function on $G_m$. 
Since $d_{\mathbf{G}_m}(g^{-1} h, h)=d_{\mathbf{G}_m}(e_{G_m}, g)$,
we have that
\begin{itemize}
\item
$\Phi_m(g^{-1} h, h)
\le (2 - \tilde{\rho}(d_{\mathbf{G}_m}(e_{G_m}, g))^2)/2$,
\item
if $d_{\mathbf{G}_m}(e_{G_m}, g) \le C$, then $\Phi_m(g^{-1} h, h) 
> 1 - \epsilon
$.
\end{itemize}
Therefore the positive definite function
$\phi_m$ satisfies the following conditions:
\begin{itemize}
\item
$\phi_m(e_{G_m}) = 1$,
\item
$\phi_m(g) \le (2 - \tilde{\rho}(d_{\mathbf{G}_m}(e_{G_m}, g))^2)/2$,
\item
if $d_{\mathbf{G}_m}(e_{G_m}, g) \le C$, then $\phi_m(g) 
> 1 - \epsilon
$.
\end{itemize}
Define positive definite functions $\psi_m$ on $F_k$
by $\psi_m = \phi_m \circ \varphi_m$, where $\varphi_m$
is the marked group quotient map $\mathbf{F_k} \rightarrow \mathbf{G}_m$.
By replacing with a subsequence, we may further assume that $\psi_m(f)$ converges for every $f \in F_k$.
Denote by $\psi$ the limit of $(\psi_m)_m$.

Let $\varphi_{\infty}$ be the marked group quotient map $\mathbf{F_k} \rightarrow \mathbf{G}_{\infty}$.
For every $f \in \mathrm{Ker}(\varphi_{\infty})$, there exists $m_f\in \mathbb{N}$ such that  $f\in \mathrm{Ker}(\varphi_m)$ for all $m\geq m_f$.
Hence
$\psi(f) = \lim_{m \rightarrow \infty} \psi_m(f)
= 1$.
By Lemma \ref{LemmaPushForward},
there exists a positive definite function $\phi$ on $G_{\infty}$ such that $\psi = \phi \circ \varphi_{\infty}$.
The function $\phi$ satisfies
\begin{itemize}
\item
$\phi(e_{G_{\infty}}) = 1$,
\item
$\phi(g) 
\le (2 - \tilde{\rho}(d_{\mathbf{G}_{\infty}}(e_{G_{\infty}},g))^2)/2$,
\item
if $d_{\mathbf{G}_{\infty}}(e_{G_{\infty}},g) \le C$, then 
$\phi(g) > 1 - \epsilon$,
\end{itemize}
By the second condition, the function $\phi_\infty$
is an element of $c_0(G)$.
It follows that $G$ is a-T-menable.
\end{proof}

\begin{proof}[Proof of $(ii)$ of Theorem~$\ref{mthm=theoremA}$]
If all $\mathbf{G}_m$ are finite, then it is immediate from Theorem~\ref{theorem=a-T-menable}. For the general case of amenable groups, we replace $\phi_m$ by employing an invariant mean on $G_m$; this modification is quite similar to the argument in the proof of [$(1)\Rightarrow (3)$] of Theorem~\ref{theorem=Item(i)}, and hence we leave the details to the reader.
\end{proof}

\begin{remark}
This proof uses specialty of Hilbert spaces (the notion of positive definite functions). However, by adopting a well-known \textit{trick of Gromov} (see \cite[Proposition~4.4]{deCornulierTesseraValette} and \cite[Section~9]{NaorPeres}), we are able to extend Theorem~\ref{theorem=a-T-menable} to more general cases. We will discuss that in our Part II paper \cite{MSPartII} in the context of \textit{fibred coarse embeddings}; see \cite[Theorem~A]{MSPartII} for the main theorem.
\end{remark}


\section{Examples}
\label{section=Examples}
In this section, we will prove Corollary~\ref{mcor=ChoiceOfGenerators} and Proposition~\ref{proposition=ChoiceOfGenerators}. Proofs are done by determining the Cayley limit group in each setting. 

\subsection{\textit{Local} point of view}\label{subsection=Local}
Before proceeding to the proofs, we describe the basic idea to identify the Cayley limit group by giving an (intuitive and sketchy) proof of the following example, which we learned from R. Tessera.

\begin{lemma}\label{lemma=Tessera}
The sequence $((\mathbb{Z}/(2^{m+1}+1)\mathbb{Z};[1]_{2^{m+1}+1},[m]_{2^{m+1}+1}))_{m\in \mathbb{N}_{\geq 2}}$ converges to $(\mathbb{Z}^2;(1,0),(0,1))$ in $\mathcal{G}(2)$.
\end{lemma}

\begin{proof}(Intuitive argument)
Recall that the Cayley convergence is the \textit{local} convergence among rooted diagrams; see $(\star)$ and Lemma~\ref{lemma=Neighborhood}.

Identify $\mathbb{Z}/(2^m+1)\mathbb{Z}$ with $[-2^m.2^m]$; recall our notation from Introduction. Since $(\mathbb{Z}/(2^{m+1}+1)\mathbb{Z};[1]_{2^{m+1}+1},[m]_{2^{m+1}+1})$ is isomorphic to $([-2^m.2^m];[1]_{2^{m+1}+1},[m]_{2^{m+1}+1})$, these two give the same point in $\mathcal{G}(2)$. Hence, we discuss the convergence of $((G_m;s_1^{(m)},s_2^{(m)})=([-2^m.2^m];[1]_{2^{m+1}+1},[m]_{2^{m+1}+1}))_{m\in \mathbb{N}_{\geq 2}}$. 

One key here is the following: Globally, the edge colorings in color $1$ (corresponding to $s^{(m)}$) on $\mathrm{CayD}(\mathbf{G}_m)$ may be seen as a shift on the cycle of length $2^m+1$; however, from the \textit{local} viewpoint, \textit{we cannot distinguish this shift from the shift on an infinite line}. More precisely, for $l\in \mathbb{Z}$ such that $|l|$ is small relative to $m$, from the root $[0]_{2^{m+1}+1} \in [-2^m.2^m]$, following ongoing edges in color $1$ (corresponding to $s_1^{(m)}$) for $l$ times (if $l<0$, then it means that we follow incoming edges in color $1$ in the orientation-reversing way for $-l$ times) \textit{gives exactly the same functions as we do it on an infinite line}. If we do similar moves for more than $2^m$ times, then we will realize that $[2^m+1]_{2^{m+1}+1}=[-2^m]_{2^{m+1}+1}$ and that we are not on the line. However, \textit{in the local point of view, we will never come up with them}. In terms of Cayley diagrams, the argument above corresponds to saying that in $R$-ball centered at $[0]$ of $\mathrm{CayD}(\mathbf{G}_m)$ with $R\leq R_1=2^m$, the edges colored in $1$ exactly look like them for $R$-balls of the Cayley diagram of $(\mathbb{Z};1)$. Similarly, on the edges colored in $2$ as long as $R\leq R_2=2^m/m$.

In order to consider mutual relation between the two generators $s_1^{(m)}$ and $s_2^{(m)}$, we proceed our local viewpoints further. Observe that $(s_1^{(m)})^m=s_2^{(m)}$; this amounts to saying that in the Cayley diagram, following the outgoing edge colored in $1$ for $m$ times results in the same function as following the outgoing edge colored in $2$. However, if we restrict ourselves to balls of radius no more than $R_3=\lfloor(m-1)/2 \rfloor$, then we will \textit{not} observe it. On the other hand, \textit{the following will be observed even in a local picture}: $s_1^{(m)}s_2^{(m)}=s_2^{(m)}s_1^{(m)}$ holds. In the Cayley diagrams, it correspond to saying that following the outgoing edge colored in $1$ and then following one colored in $2$ results in the same function as following the outgoing edge colored in $2$ and then following one colored in $1$. This will be observed as soon as our radius of ball is at least $2$. To summarize, from the \textit{local} point of view, the function of $s_1^{(m)}$ in the Cayley diagram(edges in color $1$) is \textit{indistinguishable to a shift on an infinite line}, and the same holds for the function of $s_2^{(m)}$; their functions \textit{commute}, but there is \textit{no other relation on them}.

From our discussions above, we conclude that in local point of view (namely, if $R=R_m$ is sufficiently small than $m$ but with $\lim_{m\to \infty}R_m=\infty$), $B_{\mathrm{CayD}(\mathbf{G}_m)}([0]_{2^{m+1}+1},R)$ is isomorphic, as a rooted diagram, to $B_{\mathrm{CayD}(\mathbb{Z}^2;(1,0),(0,1))}(0,R)$. By $(\star)$, this exactly means that $\mathbf{G}_m$ converges to $(\mathbb{Z}^2;(1,0),(0,1))$.

The best possible $R=R_m$ above is $\min\{R_1,R_2,R_3\}=\lfloor(m-1)/2 \rfloor$. This $R_m \to\infty$ as $m\to \infty$. However, since the diameter of $\mathrm{CayD}(\mathbf{G}_m)$ is $\asymp 2^{m+1}/m$, \textit{our $R_m$-ball is quite small $($``local''$)$ compared with the whole $($``global''$)$ picture of $\mathrm{CayD}(\mathbf{G}_m)$}.
\end{proof}

In the discussions in proceeding examples, we will only describe how local viewpoints are employed in an intuitive way to identify the Cayley limit group. However, everything in our proof may be rewritten in a completely rigorous manner.

\subsection{Symmetric groups}\label{subsection=SymmetricGroup}
We study the example $((\mathrm{Sym}(m);\sigma^{(m)},\tau^{(m)}))_{m\in \mathbb{N}_{\geq 3}}$ as in Example~\ref{example=SymmetricGroup}. We may assume that $m$ is odd, and write $m=2n+1$. Then as a marked group
\[
(\mathrm{Sym}(m);\sigma^{(m)},\tau^{(m)})\cong(\mathrm{Sym}([-n.n]);(01),(-n \cdots n)),
\]
where the second generator in the right-hand side above is the cycle permutation. Therefore, we consider the convergence of marked groups as in the right-hand side of this isomorphism, and we reset as $\sigma^{(m)}=(01)$ and $\tau^{(m)}=(-n\cdots n)$. (If $m$ is even, we use the isomorphism $\mathrm{Sym}(m)\simeq \mathrm{Sym}([-n+1.n])$, where $n=m/2$.) Then from a \textit{local} point of view, the function of $\tau^{(m)}$ in the Cayley diagram is \textit{indistinguishable to that of a shift $($infinite cyclic permutation$)$ on the infinite line $\mathbb{Z}$}; $\sigma^{(m)}$ may be seen as the transposition between $0$ and $1$ on $\mathbb{Z}$. From this reasoning, define as
\[
\sigma^{(\infty)}=(01),\ \tau^{(\infty)}=\textrm{$($the shift on $\mathbb{Z}$ by $+1$$)$} \quad \in \mathrm{Sym}(\mathbb{Z}).
\]
Then it is easy to show that 
\[
\langle \sigma^{(\infty)},\tau^{(\infty)}\rangle =\mathrm{Sym}_{<\aleph_0}(\mathbb{Z})\rtimes \mathbb{Z},
\]
where $\tau^{(\infty)}$ is the generator of $\mathbb{Z}$ in the semidirect product above; recall our notation of $\mathrm{Sym}_{<\aleph_0}(B)$ from Introduction. Therefore, we have the following:

\begin{proposition}\label{proposition=LimitSymmetricGroup}
For $((\mathrm{Sym}(m);\sigma^{(m)},\tau^{(m)}))_{m\in \mathbb{N}_{\geq 3}}$ as in Example~$\ref{example=SymmetricGroup}$, 
\begin{eqnarray*}
\lim_{m\to \infty}(\mathrm{Sym}(m);\sigma^{(m)},\tau^{(m)})=(\mathrm{Sym}_{<\aleph_0}(\mathbb{Z})\rtimes \mathbb{Z};\sigma^{(\infty)},\tau^{(\infty)}).
\end{eqnarray*}
\end{proposition}

\begin{remark}\label{remark=SymmetricGroup}
This construction may be extended to LEA groups as follows. Suppose that $(\mathbf{G}_m=(G_m;t_1^{(m)},\ldots ,t_k^{(m)}))_{m\in \mathbb{N}}$ is a LEA approximation of $\mathbf{G}_{\infty}=(G_{\infty};t_1^{(\infty)},\ldots ,t_k^{(\infty)})$. Without loss of generality, we may assume that $t_j^{(m)}\ne e_{G_m}$ for every $m\in \mathbb{N}\cup \{\infty\}$ and every $j\in [k]$. For a countable group $G$ and for $\gamma \in G\setminus \{e_G\}$, set elements $\chi_{\gamma}$ and $\theta_{\gamma}$ in $\mathrm{Sym}(G_m)$ as 

\begin{eqnarray*}
\chi_{\gamma}&=&(\textrm{the transposition on $\{e_{G},\gamma\}$}), \\
\theta_{\gamma}&=&(\textrm{the permutation on $G$ given by the right-multiplication of $\gamma$}).
\end{eqnarray*}
It is then not difficult to show that for $m\in \mathbb{N}\cup\{\infty\}$,
\begin{eqnarray*}
\Theta_m=\langle \chi_{t_1^{(m)}}, \ldots ,\chi_{t_k^{(m)}},\theta_{t_1^{(m)}},\ldots ,\theta_{t_k^{(m)}}\rangle \quad \textrm{equals}\quad
\left\{
\begin{array}{cll}
\mathrm{Sym}(G_m), & \textrm{if}\quad \sharp(G_m)<\infty,\\
\mathrm{Sym}_{<\aleph_0}(G_m)\rtimes G_m, 
& \textrm{if}\quad \sharp(G_m)=\infty.
\end{array}
\right.
\end{eqnarray*}
Indeed, observe that for every $\gamma\in G_m\setminus \{e_{G_m}\}$, the element $\chi_{\gamma}\in \mathrm{Sym}_{<\aleph_0}(G_m)$ may be written as some product of the marking above.

Thus, we obtain the following Cayley convergent sequence of groups in $\mathcal{G}(2k)$; see \cite[Lemma~4.10]{MimuraRF} for the proof.
\begin{eqnarray*}
(\Theta_m;\chi_{t_1^{(m)}}, \ldots ,\chi_{t_k^{(m)}},\theta_{t_1^{(m)}},\ldots ,\theta_{t_k^{(m)}}) &\stackrel{\mathrm{Cay}}{\longrightarrow}& (\Theta_{\infty};\chi_{t_1^{(\infty)}}, \ldots ,\chi_{t_k^{(\infty)}},\theta_{t_1^{(\infty)}},\ldots ,\theta_{t_k^{(\infty)}})
\end{eqnarray*}
Since $\Theta_m$ is amenable for all $m\in \mathbb{N}$, Theorem~\ref{mthm=theoremA} applies to this sequence.
\end{remark}

These limit groups are \textit{not} RF because it contains an isomorphic copy of $\mathrm{Alt}(\mathbb{Z})$.

\begin{remark}\label{remark=AlternatingGroup}
We may construct a similar example to Example~\ref{example=SymmetricGroup} in the context of alternating groups as follows: For \textit{odd} $m\in \mathbb{N}_{\geq 3}$, let $U_m=((123),(123\cdots m))$. Since $m$ is odd, this is a $2$-marking of $\mathrm{Alt}(m)$, the alternating group on $[m]$. Then, a similar argument as above shows that as odd $m\to \infty$, 
\begin{eqnarray*}
(\mathrm{Alt}(m);U_m) &\stackrel{\mathrm{Cay}}{\longrightarrow}& (\mathrm{Alt}(\mathbb{Z}) \rtimes \mathbb{Z};U_{\infty}),
\end{eqnarray*}
for some ``limit marking'' $U_{\infty}$. The Cayley limit group $\mathrm{Alt}(\mathbb{Z}) \rtimes \mathbb{Z}$ is amenable. 

We also note that, by employing the Nielsen--Schreier algorithm on generators of finite index subgroups, we may have a Cayley convergent sequence of alternating groups out of that of symmetric groups; see \cite[Lemma~4.11]{MimuraRF} for more details.

\end{remark}

\subsection{Unimodular and special linear groups}
\label{subsection=ExtendedLinearGroups}
Let $A$ be a commutative, associative ring with $1$. For a non-empty set $B$, consider the semigroup of all matrices $(a_{i,j})_{i,j\in B}$ over $A$ such that $a_{i,j}=0$ all but finitely many $j$ for every fixed $i$ and that $a_{i,j}=0$ all but finitely many $i$ for every fixed $j$; this is in fact a monoid with unit $I$ (the identity matrix). We define $\mathrm{GL}(B,A)$ to be the group of all invertible elements of this monoid. We warn that if $\sharp (B)=\infty$, then this group is larger than $\mathrm{GL}_{<\aleph_0}(B,A)$, which is in this paper defined as the union of all $\mathrm{GL}(K,A)$ over all non-empty finite subsets $K\subseteq B$ through the natural inclusion $\mathrm{GL}(K,A)\hookrightarrow \mathrm{GL}(B,A)$; in some literature, the symbol ``$\mathrm{GL}(B,A)$'' is used for the latter group $\mathrm{GL}_{<\aleph_0}(B,A)$. We define the  subgroup $\mathrm{SL}^{\pm}(B,A)$, the \textit{unimodular linear group}, to be the union of all $\mathrm{SL}^{\pm}(K,A)$ over all non-empty finite subsets $K\subseteq B$ (through $\mathrm{GL}(K,A)\hookrightarrow \mathrm{GL}(B,A)$). Same as what we mentioned in Introduction, for a non-empty finite subset $K\subseteq B$, define
\[
\mathrm{SL}^{\pm}(K,A)=\{g\in \mathrm{GL}(K,A): \mathrm{det}(g)\in \{\pm 1\}\}.
\]
If $B$ has a total order $\succ$, then we can define another subgroup $N_{\succ}^{\pm}(B,A)$ of $\mathrm{GL}(B,A)$ as follows: For each non-empty finite subset $K\subseteq B$, define 
\[
N_{\succ}^{\pm}(K,A)=\{(a_{i,j})_{i,j\in K}: a_{i,i}\in \{\pm 1\} \textrm{ for all $i\in K$},\ 
a_{i,j}=0 \textrm{ for all $i,j\in K$ with $i\succ j$}\}.
\]
(This is the group of ``upper triangular matrices with $\pm 1$ on diagonals.'') We define $N_{\succ}^{\pm}(B,A)$ to be the union of such $N_{\succ}^{\pm}(K,A)$ via $\mathrm{GL}(K,A)\hookrightarrow \mathrm{GL}(B,A)$. Note that each $N_{\succ}^{\pm}(K,A)$ is virtually nilpotent for such $K$; therefore $N_{\succ}^{\pm}(B,A)$ is amenable. For $i,j \in B$ with $i\ne j$ and $a\in A$, we define an \textit{elementary matrix} $e_{i,j}^a$ by
\begin{eqnarray*}
(e_{i,j}^a)_{k,l}=
\left\{
\begin{array}{cll}
1, & \quad \textrm{for}\quad k=l,\\
a, & \quad \textrm{for}\quad (k,l)=(i,j),\\
0, & \quad \textrm{otherwise}.
\end{array}
\right.
\end{eqnarray*}
This is an element in $\mathrm{SL}^{\pm}(B,A)$. Note the following commutator calculus:
\[
[e_{i,j}^{a_1},e_{j,k}^{a_2}]=e_{i,k}^{a_1a_2} \textrm{\quad for \textit{distinct} $i,j,k\in B$ and for $a_1,a_2\in A$,} \tag{\#}
\]
where $[g,h]=ghg^{-1}h^{-1}$.

\begin{remark}\label{remark=Generators}
We provide a proof of the fact that in Example~\ref{example=SpecialLinear}, $S_m$ and $T_m$ both generate $G_m$; the proof of the corresponding statement in Example~\ref{example=Point(c)} goes along the same line. Since $S_m\subseteq T_m$, it suffices to show that $S_m$ generates $G_m$.

Observe that by taking conjugations of  $\sigma^{(m)}$ by powers of $\tau^{(m)}$, we obtain all elements of the form $e_{i,i+1}^1$ in $G_m$, where we regard $e_{m,m+1}^1$ as $e_{m,1}^1$. By $(\#)$, it then follows that we may express every element of the form $e_{i,j}^r$, where $i\ne j\in [m]$ and $r\in \mathbb{F}_{p^{n_m}}$, as a certain product of $\sigma^{(m)}$, $\upsilon^{(m)}$ and $\tau^{(m)}$. Indeed, $\{1,t_{n_m}\}$ generates $\mathbb{F}_{p^{n_m}}$ as a ring. The Gaussian elimination process implies that these elements of the form above altogether generate $\mathrm{SL}(m,\mathbb{F}_{p^{n_m}})$; finally, $\mathrm{SL}(m,\mathbb{F}_{p^{n_m}})$ together with $\delta^{(m)}$ generates the unimodular group $G_m$.
\end{remark}

First we discuss the Cayley limit groups of two sequences in Example~\ref{example=Point(c)}. In a similar argument to one in Subsection~\ref{subsection=SymmetricGroup}, we observe the following: Define elements in $\mathrm{GL}(\mathbb{Z},\mathbb{Z})$ as 
\begin{eqnarray*}
&\sigma^{(\infty)}&=e_{0,1}^1,\quad \sigma'^{(\infty)}=e_{1,0}^1, \\
&\delta^{(\infty)}&=\textrm{$($diagonal matrix $-1$ on $(0,0)$-entry and $1$ on the other diagonals$)$},\\
&\tau^{(\infty)}&=\textrm{$($permutation matrix associated to the shift on $\mathbb{Z}$ by $+1$$)$},
\end{eqnarray*}
and define two subgroups $\Theta_{\infty},\Lambda_{\infty} \leqslant \mathrm{GL}(\mathbb{Z},\mathbb{Z})$ by
\[
\Theta_{\infty}=\langle \sigma^{(\infty)},\delta^{(\infty)},\tau^{(\infty)}\rangle \quad \textrm{and} \quad \Lambda_{\infty}=\langle \sigma^{(\infty)},\sigma'^{(\infty)},\delta^{(\infty)},\tau^{(\infty)}\rangle .
\]
Then, it follows that
\begin{eqnarray*}
(\mathrm{SL}^{\pm}(m,\mathbb{Z}/l_m\mathbb{Z});\sigma^{(m)},\delta^{(m)},\tau^{(m)}) &\stackrel{\mathrm{Cay}}{\longrightarrow}& (\Theta_{\infty};\sigma^{(\infty)},\delta^{(\infty)},\tau^{(\infty)}), \\
(\mathrm{SL}^{\pm}(m,\mathbb{Z}/l_m\mathbb{Z});\sigma^{(m)},\sigma'^{(m)},\delta^{(m)},\tau^{(m)}) &\stackrel{\mathrm{Cay}}{\longrightarrow}& (\Lambda_{\infty};\sigma^{(\infty)},\sigma'^{(m)},\delta^{(\infty)},\tau^{(\infty)}).
\end{eqnarray*}

\begin{proposition}\label{proposition=LimitLinear}
In the setting above, 
\begin{itemize}
 \item The group $\Theta_{\infty}$ equals $N_{>}^{\pm}(\mathbb{Z},\mathbb{Z})\rtimes  \mathbb{Z}$, where $\mathbb{Z}$ in the semidirect product acts by infinite shifts $($generated by the shift by $+1$$)$ on the coordinate set $\mathbb{Z}$. Here $>$ is the standard total order on $\mathbb{Z}$. This group is $\mathrm{amenable}$. 
 \item The group $\Lambda_{\infty}$ equals $\mathrm{SL}^{\pm}(\mathbb{Z},\mathbb{Z})\rtimes \mathbb{Z}$, where $\mathbb{Z}$ in the semidirect product acts by infinite shifts $($generated by the shift by $+1$$)$ on $\mathbb{Z}$. This group is $\mathrm{not}$ a-$\mathrm{T}$-menable.
\end{itemize}
\end{proposition}

\begin{proof}
First we verify the latter item. Since the ring $\mathbb{Z}$ is euclidean, $\mathrm{SL}(K,\mathbb{Z})$ is generated by elementary matrices over $K$ for every non-empty finite subset $K\subseteq \mathbb{Z}$ (Gaussian elimination). This proves that $\Lambda_{\infty}=\mathrm{SL}^{\pm}(\mathbb{Z},\mathbb{Z})\rtimes \mathbb{Z}$; see also $(\#)$. Because this group contains an isomorphic copy of $\mathrm{SL}(3,\mathbb{Z})$, which is well known to have \textit{Kazhdan's property $(\mathrm{T})$}, $\Lambda_{\infty}$ is not a-$\mathrm{T}$-menable; see \cite{BookBekkadelaHarpeValette} and \cite{CCJJV}.

To prove the former item, observe that all conjugations of $\sigma^{(\infty)}$ and $\delta^{(\infty)}$ by powers of $\tau^{(\infty)}$ lie in $N_{>}^{\pm}(\mathbb{Z},\mathbb{Z})$; this is because the infinite shift by $+1$ on $\mathbb{Z}$ preserves the total order $>$. We can furthermore show that $\Theta_{\infty}=N_{>}^{\pm}(\mathbb{Z},\mathbb{Z})\rtimes  \mathbb{Z}$. Since we argued above that $N_{>}^{\pm}(\mathbb{Z},\mathbb{Z})$ is amenable, so is $\Theta_{\infty}$.
\end{proof}

\begin{proof}[Proof of Proposition~$\ref{proposition=ChoiceOfGenerators}$]
Apply Theorem~\ref{mthm=theoremA} to Proposition~\ref{proposition=LimitLinear}.
\end{proof}

The most essential part of the construction as in Proposition~$\ref{proposition=ChoiceOfGenerators}$ is the following: For $m\in \mathbb{N}_{\geq 3}$, the permutation by $\tau^{(m)}$ does not ``preserve the order (upper/lower)'' \textit{from a global point of view}, because this is a shift on a finite cycle. That is the reason why $\sigma^{(m)}$, $\delta^{(m)}$ and $\tau^{(m)}$ generate the whole $\mathrm{SL}^{\pm}$ when $m\in \mathbb{N}_{\geq 3}$. However, \textit{from a local point of view, this shift is indistinguishable from the shift on an infinite line} $\mathbb{Z}$, which  \textit{preserves the standard order}. Thus, when $m=\infty$, $\sigma^{(\infty)}$, $\delta^{(\infty)}$ and $\tau^{(\infty)}$ \textit{only generate a small $($amenable$)$ subgroup} of $\mathrm{SL}^{\pm}(\mathbb{Z},\mathbb{Z})\rtimes \mathbb{Z}$, ``the \textit{upper} triangular matrices'' semidirected with $\mathbb{Z}$. If we replace $(\mathbb{Z}/l_m\mathbb{Z})_{m}$ with another sequence of rings $(A_n)_n$, then the Cayley boundary will change according to the ``accumulation points of marked rings''. Nevertheless, the basic structure of limit groups is quite similar to that in Proposition~\ref{proposition=LimitLinear}. See the proof of Corollary~\ref{mcor=ChoiceOfGenerators} for the case of Example~\ref{example=SpecialLinear}.

\begin{remark}\label{remark=deLaatMimuradelaSalle}
A certain large class $\mathcal{E}$ of Banach spaces (``$(1.1)$ for some $\beta<1/2$'' as in \cite{deLaatMimuradelaSalle}) was studied in a book of N. Tomczak-Jaegermann \cite{BookTomczak-Jaegermann} and a paper of T. de Laat and M. de la Salle \cite{deLaatdelaSalle}. For the definition, see also \cite[Example~4.11.(5)]{MSPartII}.
For this class $\mathcal{E}$, in \cite{deLaatMimuradelaSalle}, de Laat--Mimura--de la Salle showed the following: For each $E\in \mathcal{E}$, there exists $n_E\in \mathbb{N}_{\geq 3}$ such that $\mathrm{SL}(n,\mathbb{Z})$ has the fixed point property with respect to all affine isometric actions on $E$. Since $\Lambda_{\infty}$ contains all copies of $\mathrm{SL}(n,\mathbb{Z})$, $n\in \mathbb{N}_{\geq 3}$, we obtain non-((fibred) coarse) embeddability results of $W$ as in Proposition~\ref{proposition=ChoiceOfGenerators} into each $E$ in that class $\mathcal{E}$; see also Remark~\ref{remark=BetaLessThan1/2}. We will discuss it in the Part II paper \cite[Corollary~B and Corollary~2.2]{MSPartII} in details; see also \cite[Subsection~9.1]{MSPartII}.
\end{remark}

\begin{remark}\label{remark=LinearGroup}
Similar to Remark~\ref{remark=SymmetricGroup}, we can generalize this framework to a LEA approximation, at least for the latter sequence. For simplicity, here we only discuss a LEF approximation. Take $(\mathbf{G}_m=(G_m;t_1^{(m)},\ldots ,t_k^{(m)}))_{m\in \mathbb{N}}$ a LEF approximation of an infinite marked group $\mathbf{G}_{\infty}=(G_{\infty};t_1^{(\infty)},\ldots ,t_k^{(\infty)})$. Without loss of generality, we may assume that $t_j^{(m)}\ne e_{G_m}$ for $m \in \mathbb{N}\cup\{\infty\}$ and every $j\in [k]$. For a countable group $G$, for $\gamma\in G\setminus \{g\}$ and for $l\in \mathbb{N}$, define elements $\sigma_{\gamma}=\sigma_{\gamma}(l)$, $\sigma'_{\gamma}=\sigma'_{\gamma}(l)$, $\delta=\delta(l)$ and $\tau_{\gamma}=\tau_{\gamma}(l)$ in $\mathrm{GL}(G,\mathbb{Z}/l\mathbb{Z})$ by 
\begin{eqnarray*}
&\sigma_{\gamma}&=e_{e_{G},\gamma}^1,\quad \sigma'_{\gamma}=e_{\gamma,e_{G}}^1, \\
&\delta&=\textrm{$($diagonal matrix $-1$ on $(e_{G},e_{G})$-entry and $1$ on the other diagonals$)$},\\
&\tau_{\gamma}&=\textrm{$($the permutation matrix by the shift on $G$ by the right multiplication of $\gamma$$)$}.
\end{eqnarray*}
More precisely, the element $\tau_{\gamma}$ is defined by
\[
(\tau_{\gamma})_{g,h}=\left\{
\begin{array}{cl}
 1, & \quad \textrm{if $g=h\gamma$},\\
 0, & \quad \textrm{otherwise},
\end{array}\right. \quad \textrm{for every $g,h\in G$.}
\]
Let $(l_m)_{m\in \mathbb{N}}$ be a sequence of integers $\in \mathbb{N}_{\geq 2}$ that satisfies $\lim_{m\to \infty}l_m=\infty$ and let $l_{\infty}=0$. 
For $m\in \mathbb{N}\cup\{\infty\}$, let $\Lambda_m$ be the group generated by $(\sigma_{t_j^{(m)}}(l_m))_{j\in [k]}$, $(\sigma'_{t_j^{(m)}}(l_m))_{j\in [k]}$, $\delta(l_m)$ and $(\tau_{t_j^{(m)}}(l_m))_{j\in [k]}$. Then for all $m\in \mathbb{N}\cup\{\infty\}$,
\begin{eqnarray*}
\Lambda_m=
\left\{
\begin{array}{cll}
\mathrm{SL}^{\pm}(G_m,\mathbb{Z}/l_m\mathbb{Z}), & \textrm{if}\quad m\in \mathbb{N},\\
\mathrm{SL}^{\pm}(G_{\infty},\mathbb{Z})\rtimes G_{\infty}, 
& \textrm{if}\quad m=\infty ,
\end{array}
\right.
\end{eqnarray*}
and $\Lambda_m$ converges to $\Lambda_{\infty}$ in $\mathcal{G}(3k+1)$ with respect to the corresponding markings. 

Now, assume that $G_{\infty}$ is \textit{right-orderable}, namely, there exists a total order $\succ$ on $G_{\infty}$ that is preserved by the right multiplications of $G$. (A group is right-orderable if and only if it is \textit{left-orderable}; the latter terminology is more common in the literature.) We refer to a monograph \cite{bookDeroinNavasRivas} for details on this concept. Fix such a (right-invariant) order $\succ$ on $G_{\infty}$. We moreover assume that $t_j^{(\infty)}\succ e_{G_{\infty}}$ for all $j\in [k]$; if this is not the case, then for each $j\in [k]$ such that $t_j^{(\infty)}\prec e_{G_{\infty}}$, replace $t_j^{(m)}$ with $(t_j^{(m)})^{-1}$ for all $m\in \mathbb{N}\cup \{\infty\}$. For $m\in \mathbb{N}\cup\{\infty\}$, let $\Theta_m$ be the group generated by $(\sigma_{t_j^{(m)}})_{j\in [k]}$, $\delta$ and $(\tau_{t_j^{(m)}})_{j\in [k]}$. Then $\Theta_m$ converges to $\Theta_{\infty}$ in $\mathcal{G}(2k+1)$ with respect to the corresponding markings. For all $m\in \mathbb{N}$,  $\Theta_m$ equals $\mathrm{SL}^{\pm}(G_m,\mathbb{Z}/l_m\mathbb{Z})$. For $m=\infty$, $\Theta_{\infty}$ is a (in general, possibly proper) subgroup of $N^{\pm}_{\succ}(G_{\infty},\mathbb{Z})\rtimes G_{\infty}$ and $\Theta_{\infty}$ contains a copy of $G_{\infty}$. In particular, $\Theta_{\infty}$ is amenable if and only if $G_{\infty}$ is amenable.

\end{remark}

Finally, we deal with the case as in Example~\ref{example=SpecialLinear}.

\begin{proof}[Proof of Corollary~$\ref{mcor=ChoiceOfGenerators}$]
Item $(i)$ follows from Proposition~\ref{proposition=LimitSymmetricGroup} and Theorem~\ref{mthm=theoremA}. 

We then discuss items $(ii)$ and $(iii)$. Let $G_m=\mathrm{SL}^{\pm}(m,\mathbb{F}_{p^{n_m}}),S_m, T_m$ be as in Example~\ref{example=SpecialLinear}. In a manner quite similar to one as in the proof of Proposition~\ref{proposition=LimitLinear}, we can prove that 
\begin{eqnarray*}
(\mathrm{SL}^{\pm}(m,\mathbb{F}_{p^{n_m}});S_m) &\stackrel{\mathrm{Cay}}{\longrightarrow}& (N_{>}^{\pm}(\mathbb{Z},\mathbb{F}_p[t])\rtimes \mathbb{Z};\sigma^{(\infty)},\upsilon^{(\infty)},\delta^{(\infty)},\tau^{(\infty)}), \\
(\mathrm{SL}^{\pm}(m,\mathbb{F}_{p^{n_m}});T_m) &\stackrel{\mathrm{Cay}}{\longrightarrow}& (\mathrm{SL}^{\pm}(\mathbb{Z},\mathbb{F}_p[t])\rtimes \mathbb{Z};\sigma^{(\infty)},\sigma'^{(\infty)},\upsilon^{(\infty)},\delta^{(\infty)},\tau^{(\infty)}), 
\end{eqnarray*}
for certain ``limit markings.'' Here $\mathbb{F}_p[t]$ denotes the one-variable polynomial ring over $\mathbb{F}_p$, and $>$ is the standard total order on $\mathbb{Z}$. More precisely, to see that the ``limit coefficient ring'' is $\mathbb{F}_p[t]$, observe that all of $\mathbb{F}_{p^{n_m}}$ are quotient rings of $\mathbb{F}_p[t]$ by $t\mapsto t_{n_m}$, and that all proper quotient rings of $\mathbb{F}_p[t]$ are finite. 

It is well known that $\mathrm{SL}(3,\mathbb{F}_p[t])$ has property $(\mathrm{T})$; this group is a lattice in $\mathrm{SL}(3,\mathbb{F}_p((t^{-1})))$, where $\mathbb{F}_p((t^{-1}))$ denotes the (local) field of formal Laurent series  with indeterminate  $t^{-1}$ over $\mathbb{F}_p$. Therefore, a similar reasoning to one in the proof of Proposition~\ref{proposition=LimitLinear}, together with Theorem~\ref{mthm=theoremA} end our proof.
\end{proof}

\begin{remark}\label{remark=Lafforgue}
Much more than having property $(\mathrm{T})$ is known for $\mathrm{SL}(3,\mathbb{F}_p[t])$: Results of Lafforgue \cite{Lafforgue1}, \cite{Lafforgue2} imply (via $L_2$-induction argument) that $\mathrm{SL}(3,\mathbb{F}_p([t])$ has the fixed point properties with respect to all affine isometric actions on complex Banach spaces with non-trivial (linear) type. In the next section, we will see a byproduct of it, Theorem~\ref{mthm=PoorCompression}; see also Theorem~\ref{theorem=Lafforgue} and Theorem~\ref{theorem=Tau}.
\end{remark}

\begin{remark}\label{remark=SpecialLinearGroup}
We may construct a similar example to Examples~\ref{example=SpecialLinear} and \ref{example=Point(c)} for special linear groups with dropping $\delta^{(m)}$ from the markings. For \textit{odd} $m\in \mathbb{N}_{\geq 3}$, $(\sigma^{(m)},\upsilon^{(m)},\tau^{(m)})$ as in Examples~\ref{example=SpecialLinear} is a $3$-marking of $\mathrm{SL}(m,\mathbb{F}_{p^{n_m}})$; recall Remark~\ref{remark=Generators}. Similarly, the quadruple $(\sigma^{(m)},\sigma'^{(m)},\upsilon^{(m)}, \tau^{(m)})$ is a $4$-marking of that. We can see that as odd $m\to \infty$,
\begin{eqnarray*}
 (\mathrm{SL}(m,\mathbb{F}_{p^{n_m}});\sigma^{(m)},\upsilon^{(m)},\tau^{(m)})&\stackrel{\mathrm{Cay}}{\longrightarrow}& \mathrm{N}_{>}(\mathbb{Z},\mathbb{F}_p[t])\rtimes \mathbb{Z}, \\
 (\mathrm{SL}(m,\mathbb{F}_{p^{n_m}});\sigma^{(m)},\sigma'^{(m)},\upsilon^{(m)},\tau^{(m)})&\stackrel{\mathrm{Cay}}{\longrightarrow}& \mathrm{SL}(\mathbb{Z},\mathbb{F}_p[t])\rtimes \mathbb{Z}.
\end{eqnarray*}
Here $\mathrm{N}_{>}(\mathbb{Z},\mathbb{F}_p[t])$ denotes the set of all upper triangular matrices (in terms of the standard order $>$ on $\mathbb{Z}$) with coefficients in $\mathbb{F}_p[t]$ \textit{with $1$ on diagonal} of finite size; $\mathrm{SL}(\mathbb{Z},\mathbb{F}_p[t])$ means the union of all special linear groups over $\mathbb{F}_p[t]$ of finite size. Again, the former Cayley limit group is amenable and the latter is not a-$\mathrm{T}$-menable.

Similarly, in the setting of Example~\ref{example=Point(c)}, as \textit{odd} $m\to \infty$,
\begin{eqnarray*}
 (\mathrm{SL}(m,\mathbb{Z}/l_m\mathbb{Z}));\sigma^{(m)},\tau^{(m)})&\stackrel{\mathrm{Cay}}{\longrightarrow}& \mathrm{N}_{>}(\mathbb{Z},\mathbb{Z})\rtimes \mathbb{Z}, \\
 (\mathrm{SL}(m,\mathbb{Z}/l_m\mathbb{Z}));\sigma^{(m)},\sigma'^{(m)},\tau^{(m)})&\stackrel{\mathrm{Cay}}{\longrightarrow}& \mathrm{SL}(\mathbb{Z},\mathbb{Z})\rtimes \mathbb{Z}.
\end{eqnarray*}

Let
\[
\mathbb{N}_{\mathrm{odd}}=\{3,5,7,\ldots\}
\]
be the set of \textit{odd} integers \textit{at least} $3$. (This is a non-standard notation, but we use it for simplicity.) 
In the setting above, if we set
\begin{eqnarray*}
X'&=& \bigsqcup_{m\in \mathbb{N}_{\mathrm{odd}}}\mathrm{Cay}(\mathrm{SL}(m,\mathbb{F}_{p^{n_m}});\sigma^{(m)},\upsilon^{(m)},\tau^{(m)}),\\
Y'&=& \bigsqcup_{m\in \mathbb{N}_{\mathrm{odd}}}\mathrm{Cay}(\mathrm{SL}(m,\mathbb{F}_{p^{n_m}});\sigma^{(m)},\sigma'^{(m)},\upsilon^{(m)},\tau^{(m)}),\\
Z'&=& \bigsqcup_{m\in \mathbb{N}_{\mathrm{odd}}}\mathrm{Cay}(\mathrm{SL}(m,\mathbb{Z}/l_m\mathbb{Z});\sigma^{(m)},\tau^{(m)}),\\
W'&=& \bigsqcup_{m\in \mathbb{N}_{\mathrm{odd}}}\mathrm{Cay}(\mathrm{SL}(m,\mathbb{Z}/l_m\mathbb{Z});\sigma^{(m)},\sigma'^{(m)},\tau^{(m)}),
\end{eqnarray*}
then by Theorem~\ref{mthm=theoremA}, $X'$ and $Z'$ have property A. On the other hand, $Y'$ and $W'$ fail to have a coarse embedding into a Hilbert space; we will study more on (fibred) coarse non-embeddability of $Y'$ and $W'$ in our Part II paper \cite[Corollary~2.2 and Subsection~9.1]{MSPartII}.

If $m\in \mathbb{N}_{\geq 3}$ is even, then we may take the following modification: in this case, $\tau^{(m)}$ itself does \textit{not} lie in $\mathrm{SL}(m,\mathbb{F}_{p^{n_m}})$. Instead, $\delta^{(m)}\tau^{(m)}$ belongs to $\mathrm{SL}(m,\mathbb{F}_{p^{n_m}})$. In a similar argument to one in Remark~\ref{remark=Generators}, it follows that $(\sigma^{(m)},\upsilon^{(m)},\delta^{(m)}\tau^{(m)})$ is a $3$-marking of $\mathrm{SL}(m,\mathbb{F}_{p^{n_m}})$. Then, by Lemma~\ref{lemma=LEFsubgroups}, the sequence $((\mathrm{SL}(m,\mathbb{F}_{p^{n_m}};\sigma^{(m)},\upsilon^{(m)},\delta^{(m)}\tau^{(m)}))_{m}$ converges in the Cayley topology to the marked group with marking $(\sigma^{(\infty)},\upsilon^{(\infty)},\delta^{(\infty)}\tau^{(\infty)})$. It is easy to see that the group generated by $\{\sigma^{(\infty)},\upsilon^{(\infty)},\delta^{(\infty)}\tau^{(\infty)}\}$ equals $\mathrm{N}_{>}(\mathbb{Z},\mathbb{F}_p[t])\rtimes \mathbb{Z}$. A similar argument to one above applies to the case where we add $\sigma'{}^{(m)}$ to the marking. Therefore, we have the convergence as \textit{even} number $m\to \infty$,
\begin{eqnarray*}
 (\mathrm{SL}(m,\mathbb{F}_{p^{n_m}});\sigma^{(m)},\upsilon^{(m)},\delta^{(m)}\tau^{(m)})&\stackrel{\mathrm{Cay}}{\longrightarrow}& \mathrm{N}_{>}(\mathbb{Z},\mathbb{F}_p[t])\rtimes \mathbb{Z}, \\
 (\mathrm{SL}(m,\mathbb{F}_{p^{n_m}});\sigma^{(m)},\sigma'^{(m)},\upsilon^{(m)},\delta^{(m)}\tau^{(m)})&\stackrel{\mathrm{Cay}}{\longrightarrow}& \mathrm{SL}(\mathbb{Z},\mathbb{F}_p[t])\rtimes \mathbb{Z},
\end{eqnarray*}
with respect to certain markings of the Cayley limit groups.
\end{remark}

\begin{remark}\label{remark=PartII}
We here announce that in our Part II paper \cite[Theorem~D]{MSPartII}, we prove the following. The construction is based on Remark~\ref{remark=Lafforgue} and Remark~\ref{remark=SymmetricGroup}; see \cite[Subsection~9.6]{MSPartII}.

\begin{theorem}
There exist $(k_l)_{l\in \mathbb{N}}$ of a sequence of natural numbers at least $2$ with $\lim_{l\to \infty}k_l=\infty$ and two $($ordered$)$ systems of generators $(\Xi_l)_{l\in \mathbb{N}}$, $(\Omega_l)_{l\in \mathbb{N}}$ of symmetric groups $(\mathfrak{S}([k_l]))_{l\in \mathbb{N}}$ that satisfy all of the following.
\begin{enumerate}[$(1)$]
    \item For all $l\in \mathbb{N}$, $\sharp(\Xi_l)=8$ and $\sharp(\Omega_l)=9$. For each $l\in \mathbb{N}$, $\Omega_l$ is constructed by adding one extra element to $\Xi_l$.  
\item A coarse disjoint union $\coprod_{l\in \mathbb{N}}\mathrm{Cay}(\mathfrak{S}([k_l]);\Xi_l)$ has property A.
  \item A coarse disjoint union $\coprod_{l\in \mathbb{N}}\mathrm{Cay}(\mathfrak{S}([k_l]);\Omega_l)$ does $\mathrm{not}$ admit a fibred coarse embedding into any of these spaces: 
  \begin{itemize}
    \item Complex Banach spaces of non-trivial type, and complex Banach spaces that are sphere equivalent to Banach spaces of non-trivial type.    
\item Complete $\mathrm{CAT}(0)$ spaces with the Izeki--Nayatani invariant $\delta$ $($see Definition~$6.1$ in \cite{IzekiNayatani} for the definition$)$ being strictly smaller than $1$.
  \end{itemize}
\end{enumerate}
Moreover, this construction may be done in an explicit manner.
\end{theorem}
\end{remark}

\section{Disjoint union with property A with very poor compression function}
\label{section=PoorCompression}

We will prove Theorem~\ref{mthm=PoorCompression}. The key to the proof is the \textit{Banach spectral gap} of finite graphs. Recall our conventions of symbols $\precsim$, $\succsim$, $\asymp$ and $\precnsim$ from Introduction.

By following \cite[Definition~1.3]{MimuraSphereEquivalence}, we say that two Banach spaces $E$ and $F$ are \textit{sphere equivalent} if there exists a bijective map $\Phi\colon S(E)\to S(F)$ between unit spheres such that $\Phi$ and $\Phi^{-1}$ are both uniformly continuous. The Pisier characterization tells that a complex Banach space $E$ is \textit{of  trivial $($linear$)$ type} if and only if it admits homomorphic embeddings of $(\ell_{1,\mathbb{C}}^n)_{n\in \mathbb{N}}$ (finite dimensional complex $\ell_1$-spaces) with uniform biLipschitz constants; see \cite[Example~4.11.$(4)$]{MSPartII} for the original definition. See also \cite{BookTomczak-Jaegermann} and \cite{BookBenyaminiLindenstrauss} on comprehensive treatments of geometry of Banach spaces.

\subsection{Decay of Banach spectral gaps and distortions}
\label{SubsectionBanachGap}

\begin{definition}\label{definition=BanachSpectralGap}
Let $\Gamma=(V_{\Gamma},E_{\Gamma})$ be a finite connected (undirected) graphs (possibly with multi-edges or self-loops), where $V_{\Gamma}$ is the vertex set and $E_{\Gamma}$ is the edge set. For each $e\in E_{\Gamma}$, we put both of the possible two directions and regard it as two oriented edges. Let $\vec{E}_{\Gamma}$ be the set of all oriented edges constructed as above (hence $\sharp(\vec{E}_{\Gamma})=2\sharp(E_{\Gamma})$). For each oriented edge $\vec{e}\in \vec{E}_{\Gamma}$, let $\vec{e}_+$ and $\vec{e}_-$ be, respectively,  the end vertex and the start vertex of $\vec{e}$. We denote by $\Delta(\Gamma)$ the maximal degree of $\Gamma$, that means, the maximal number of oriented edges going out from a vertex. We view $\Gamma$ as a metric space equipped with the path metric $d_{\Gamma}$, and denote by $\mathrm{diam}(\Gamma)$ the diameter of $\Gamma$. Let $E$ be a Banach space. 
\begin{enumerate}[$(1)$]
  \item (See Definition~1.1 in \cite{MimuraSphereEquivalence}.) The (non-normalized) \textit{Banach $(E,2)$-spectral gap} of $\Gamma$, denoted by $\lambda_1 (\Gamma,E)(=\lambda_1(\Gamma,E,2))$, is defined by 
\[
\lambda_1 (\Gamma,E)=\frac{1}{2}\inf_{f\colon V_{\Gamma}\to E}\frac{\sum_{\vec{e}\in \vec{E}_{\Gamma}}\| f(\vec{e}_+)-f(\vec{e}_-)\|^2_E}{\sum_{v\in V_{\Gamma}}\| f(v)-m(f)\|^2_E}.
\]
Here $f$ moves among all non-constant maps $V\to E$, and $m(f)$ denotes the mean $(\sum_{v\in V_{\Gamma}}f(v))/\sharp(V_{\Gamma})$ of $f$.
  \item The \textit{$E$-distortion} of $\Gamma$, denoted by $c_E(G)$, is defined by 
\[
c_E(\Gamma)=\inf_{f\colon V_{\Gamma}\to E,\ \mathrm{biLipschitz}}\|f\|_{\mathrm{Lip}}\cdot\|f^{-1}\|_{\mathrm{Lip}}.
\]
Here $\|\cdot\|_{\mathrm{Lip}}$ denotes the Lipschitz constant, and $f^{-1}$ means $f^{-1}\colon f(V_{\Gamma})\to V_{\Gamma}$.
\end{enumerate}
\end{definition}

The following two results altogether connect Banach spectral gaps to compression functions of coarse embeddings.

\begin{theorem}
[Special case of the generalized Jolissaint--Valette inequality, see \cite{JolissaintValette} and Theorem~$2.6$ in \cite{MimuraSphereEquivalence}]
\label{theorem=JolissaintValette}
Let $\Gamma$ be a finite connected $\mathrm{Cayley}$ graph. Then for every Banach space $E$, it holds that
\[
c_E(\Gamma)\geq \mathrm{diam}(\Gamma) \sqrt{\frac{\lambda_1(\Gamma,E)}{2\Delta(\Gamma)}}.
\]
In particular, 
$c_E(\Gamma)\succsim_{\Delta(\Gamma)} \lambda_1(\Gamma,E)^{\frac{1}{2}} \cdot \mathrm{diam}(\Gamma)$.
\end{theorem}

\begin{lemma}
[Generalized version of Austin's Lemma (Lemma~3.1 in \cite{Austin})]
\label{lemma=Austin}
Let $(\Gamma_m)_{m\in \mathbb{N}}$ be a sequence of finite connected graphs with $\mathrm{diam}(\Gamma_m)\to \infty$. Let $\rho \colon [0,\infty) \to [0,\infty)$ be a proper function such that $\rho(r)/r$ is non-increasing for $r\in [0,\infty)$ large enough. Let $E$ be a Banach space. If for $m\in \mathbb{N}$ large enough it holds that
\[
\frac{\mathrm{diam}(\Gamma_m)}{\rho (\mathrm{diam}(\Gamma_m))} \precnsim c_E(\Gamma_m),
\]
then for $\omega(r)=r\colon [0,\infty) \to [0,\infty)$, $(\rho,\omega)\not \in \mathcal{CP}_E(\bigsqcup_{m\in \mathbb{N}}\Gamma_m)$.
\end{lemma}

For the convenience of the readers, we give a proof of Lemma~\ref{lemma=Austin}.

\begin{proof}[Proof of Lemma~$\ref{lemma=Austin}$]
Suppose, on the contrary, that there exists a coarse embedding $f\colon \bigsqcup_{m\in \mathbb{N}} \Gamma_m \to E$ for which $(\rho,\omega)$ above is a control pair. Set $f_m=f\mid_{\Gamma_m} \colon V_m \to E$ for each $m$, where $V_m$ is the vertex set of $\Gamma_m$. With bounded adjustment, we may assume that $f_m$ is injective. Since $\omega(r)=r$, for all $m$ it holds that $\|f_m\|_{\mathrm{Lip}}\leq 1$.  Then for $m$ sufficiently large, we have that
\begin{eqnarray*} 
 \frac{\mathrm{diam}(\Gamma_m)}{\rho (\mathrm{diam}(\Gamma_m))} 
&\precnsim& c_E(\Gamma_m) \\
&\leq& \|f_m^{-1}\|_{\mathrm{Lip}} \\
&\leq& \max_{v\ne w\in V_m}\frac{d(v,w)}{\|f_m(v)-f_m(w)\|} \\
&\leq& \max_{v\ne w\in V_m}\frac{d(v,w)}{\rho(d(v,w))} \\
& \precsim_{\rho}& \frac{\mathrm{diam}(\Gamma_m)}{\rho(\mathrm{diam}(\Gamma_m))},
\end{eqnarray*}
a contradiction.
\end{proof}

The goal of this section is to prove the following.

\begin{theorem}
\label{theorem=BanachSpectralGap}
Let $G_m=\mathrm{SL}^{\pm}(m,\mathbb{F}_{p^{n_m}}),\sigma^{(m)}, \upsilon^{(m)},\delta^{(m)},\tau^{(m)}$ be as in Example~$\ref{example=SpecialLinear}$. Let $\Gamma_m=\mathrm{Cay}(G_m;\sigma^{(m)}, \upsilon^{(m)},\delta^{(m)},\tau^{(m)})$
Let $E$ be a complex Banach space that is sphere equivalent to a complex Banach space of non-trivial type. Then we have that 
\[
\lambda_1(\Gamma_m,E)  \succsim_{E,p} m^{-6}.
\]
\end{theorem}
The point here is, as the symbol ''$\succsim_{E,p}$'' indicates, this estimate does \textit{not} depend on the choices of $(n_m)_{m\in \mathbb{N}_{\geq 3}}$. 

\begin{remark}\label{remark=ell_2}
In what follows, to simplify our argument, we replace a Banach space $E$ with $\ell_2(\mathbb{N},E)$ so that $E$ and $\ell_2(\mathbb{N},E)$ are isometrically isomorphic; indeed, it is easy to see that $\lambda_1(\Gamma,E)=\lambda_2(\Gamma,\ell_2(\mathbb{N},E))$ (see for instance \cite[Lemma~2.1]{MimuraSphereEquivalence}), and if a generalized metric space only admits coarse embeddings into $\ell_2(\mathbb{N},E)$ with poor compressions, then the same holds for coarse embeddings into $E$. Note that if $E$ and $F$ are sphere equivalent, then so are $\ell_2(\mathbb{N},E)$ and $\ell_2(\mathbb{N},F)$; see \cite[Proposition~3.9]{MimuraSphereEquivalence}. Furthermore, if $F$ has non-trivial type, then so does $\ell_2(\mathbb{N},F)$.
\end{remark}

\subsection{$(\tau)$-type constants and Banach spectral gaps}
We relate $\lambda_1(\Gamma,E)$ to a representation theoretic constant.

\begin{definition}
\label{definition=TauConstant}
Let $\mathbf{G}=(G;S)$ be a marked group and $E$ be a Banach space.
\begin{enumerate}[$(1)$]
  \item For an isometric linear representation $\pi \colon G\to O(E)$, where $O(E)$ deotes the group of surjective linear isometries, define $\kappa (\mathbf{G},\pi) \in [0,+\infty]$ be the supremum among all $\kappa\geq 0$ that satisfy the following: For all $\xi \in E$,
\[
\sup_{s\in S}\|\pi(s)\xi-\xi \|\geq \kappa \cdot \sup_{g\in G}\|\pi(g)\xi-\xi\|
\]
holds true.

  Note that $\kappa\leq 1$ unless $E=E^{\pi(G)}$, where $E^{\pi(G)}$ denotes the space of all $\pi(G)$-invariant vectors.
  \item The \textit{$(\tau)$-type constant associated with $(\mathbf{G},E)$} is defined as the infimum of $\kappa (\mathbf{G},\pi) $ over all isometric linear representations $\pi\colon G\to O(E)$ on $E$ \textit{that factors through a finite marked group quotient map} $\mathbf{G}\twoheadrightarrow \mathbf{H}$. In this paper, we write it as $\mathcal{K}^{(\tau)}(\mathbf{G},E)$.
\end{enumerate}
\end{definition}

There are several formulations of $(\tau)$-type constant which are mutually slightly different. In what follows, we exhibit another formulation of $(\tau)$-type constants: Suppose that $\pi\colon G\to O(E)$ factor through a marked group quotient $\varphi\colon \mathbf{G}\twoheadrightarrow \mathbf{H}$. Let $\pi_H\colon H\to O(E)$ be such that $\pi=\pi_H\circ \varphi$. Then since $H$ is finite, we can define a projection
\[
P=P_{\pi,H}\colon E\twoheadrightarrow E^{\pi(G)};\quad \xi \mapsto \frac{1}{\sharp(H)}\sum_{h\in H}\pi_H(h)\xi,
\]
which has norm $1$ (or $0$ if $E^{\pi(G)}=\{0\}$); in other words, $P=(\sum_{h\in H}\pi_H(h))/\sharp(H)$. By construction, $P$ commutes with $\pi(G)$. Thus $E$ is decomposed into 
\[
E=E^{\pi(G)}\oplus (I-P)E
\]
as $\pi(G)$-representations. In short, $E^{\pi(G)}$ is complemented as a $\pi(G)$ subrepresentation. In \cite[Definition~3.1.$(ii)$]{MimuraSphereEquivalence}, a formulation of $(\tau)$-type constants was given in the following way: $\sup_{g\in G}\|\pi(g)\xi-\xi\|$ is replaced with $\|\xi\|$ and $\xi$ runs over all $\xi \in (I-P)E$. (Also, there we only considered quasi-regular representation; however, a certain Peter--Weyl type argument shows that it will not affect the order of constants as long as $\sharp (S)$ is uniformly bounded.) This formulation is not identical to the one given in Definition~\ref{definition=TauConstant}.

The following lemma, nevertheless, shows that these two formulations are equivalent.
\begin{lemma}\label{lemma=FiniteGroupQuotient}
In the setting in the paragraph above, for every $\xi \in E$,
\[
\|\xi-P\xi\|\leq \sup_{g\in G}\|\pi(g)\xi-\xi\|\leq 2\|\xi-P\xi\|.
\]
\end{lemma}
\begin{proof}
The right inequality is directly by the decomposition $\xi=P\xi+(I-P)\xi$. The left one follows from
\[
\|\xi-P\xi\|\leq \frac{1}{\sharp(H)}\sum_{h\in H} \|\pi_H(h)\xi-\xi\|.
\]
Observe that $\sup_{g\in G}\|\pi(g)\xi-\xi\|=\sup_{h\in H}\|\pi_H(h)\xi-\xi\|$.
\end{proof}

\begin{lemma}\label{lemma=SpectralGapTau}
Let $E$ be a Banach space that is isometrically isomorphic to $\ell_2(\mathbb{N},E)$ $($compare with Remark~$\ref{remark=ell_2}$$)$. Let $\mathbf{G}=(G;S)$ be a marked group. Then for every non-trivial finite marked group quotient $\varphi\colon \mathbf{G}\twoheadrightarrow \mathbf{H}=(H;S_H)$, we have that
\[
\lambda_1(\mathrm{Cay}(\mathbf{H}),E)\geq \left(\mathcal{K}^{(\tau)}(\mathbf{G},E)\right)^2.
\]
\end{lemma}
This lemma or some variant of it might be well known to experts; nevertheless, we present the proof for the convenience of the readers.
\begin{proof}
Recall our definition of the Cayley graph from Definition~\ref{definition=CayleyDiagram}. We write down the definition of $\lambda_1(\mathrm{Cay}(\mathbf{H}),H)$:
\begin{eqnarray*}
\lambda_1 (\mathrm{Cay}(\mathbf{H}),E)&=&\frac{1}{2}\inf_{f\colon H\to E}\frac{\sum_{\vec{e}\in \vec{E}_{\mathrm{Cay}(\mathbf{H})}}\| f(\vec{e}_+)-f(\vec{e}_-)\|^2}{\sum_{h\in H}\| f(h)-m(f)\|^2}\\
&=&\inf_{f\colon H\to E}\frac{\sum_{h\in H, t\in S_H}\| f(th)-f(h)\|^2_E}{\sum_{h\in H}\| f(h)-m(f)\|^2_E}.
\end{eqnarray*}
Now we regard $f\colon H\to E$ as an element in $\ell_2(H,E)$. Consider the left-regular representation $\lambda=\lambda_{H,E}$ of $H$ on $\ell_2(H,E)$, that means, $\lambda\colon H\to O(\ell_2(H,E))$; $[\lambda_hf](x)=f(h^{-1}x)$. Set $\pi=\lambda \circ \varphi$; this is an isometric linear action of $G$ on $\ell_2(H,E)$ that factors through $\varphi$. We consider the projection $P$ associated with $\pi$; recall it from the discussion above Lemma~\ref{lemma=FiniteGroupQuotient}. In this case, $P\colon \ell_2(H,E)\twoheadrightarrow \ell_2(H,E)^{\pi(G)}$; $f\mapsto Pf$ is exactly the same as taking the mean $m(f)$ of $f\colon H\to E$, where we view $m(f)$ as a ``constant'' function in $\ell_2(H,E)$. Therefore, $\lambda_1 (\mathrm{Cay}(\mathbf{H}),E)$ is further rewritten as
\begin{eqnarray*}
\lambda_1 (\mathrm{Cay}(\mathbf{H}),E)&=&\inf_{f\in \ell_2(H,E)}\frac{\sum_{t\in S_H}\sum_{h\in H}\| f(th)-f(h)\|^2_E}{\| f-Pf\|^2_{\ell_2(H,E)}}\\
&=&\inf_{f\in \ell_2(H,E)}\frac{\sum_{t\in S_H}\| \lambda_{t^{-1}}f-f\|^2_{\ell_2(H,E)}}{\| f-Pf\|^2_{\ell_2(H,E)}}\\
&=&\inf_{f\in \ell_2(H,E)}\frac{\sum_{t\in S_H}\| \lambda_{t}f-f\|^2_{\ell_2(H,E)}}{\| f-Pf\|^2_{\ell_2(H,E)}}\\
&=&\inf_{f\in \ell_2(H,E)}\frac{\sum_{s\in S}\| \pi(s)f-f\|^2_{\ell_2(H,E)}}{\| f-Pf\|^2_{\ell_2(H,E)}}\\
&\geq& \inf_{f\in \ell_2(H,E)}\frac{\sum_{s\in S}\| \pi(s)f-f\|^2_{\ell_2(H,E)}}{\sup_{g\in G}\| \pi(g)f-f\|^2_{\ell_2(H,E)}}.
\end{eqnarray*}

Indeed, the very last inequality follows from Lemma~\ref{lemma=SpectralGapTau}. 
Here $f$ runs over all vectors in $\ell_2(H,E)$ such that $f-Pf\ne 0$; this is equivalent to saying that $\sup_{g\in G}\| \pi(g)f-f\|_{\ell_2(H,E)}\ne 0$.  Finally, apply the definition of $\mathcal{K}^{(\tau)}(\mathbf{G},E)$ to $\pi$. Here recall that $E$ and $\ell_2(H,E)$ are isometrically isomorphic by assumption.
\end{proof}

\subsection{Proof of Theorem~\ref{theorem=BanachSpectralGap}}

Our goal, Theorem~\ref{theorem=BanachSpectralGap}, will be deduced from the following. Recall from the beginning paragraph Subsection~\ref{subsection=ExtendedLinearGroups} the definition of elementary matrices $e_{i,j}^a$.

\begin{theorem}\label{theorem=Tau}
Let $p$ be a prime. For $m\in \mathbb{N}_{\geq 3}$, define a $4$-marked group $\tilde{\mathbf{G}}_m=(\tilde{G}_m;\tilde{\sigma}^{(m)},\tilde{\upsilon}^{(m)},\tilde{\delta}^{(m)},\tilde{\tau}^{(m)})$ as follows: 
\begin{itemize}
  \item $\tilde{G}_m=\mathrm{SL}^{\pm}(m,\mathbb{F}_p[t])(=\mathrm{SL}^{\pm}([m],\mathbb{F}_p[t]))$. 
  \item $\tilde{\sigma}^{(m)}=e_{1,2}^{1}$, $\tilde{\upsilon}^{(m)}=e_{1,2}^t$.
  \item $\tilde{\delta}^{(m)}$ is the diagonal matrix, $-1$ on $(1,1)$-entry and $1$ on the other diagonals.
  \item $\tilde{\tau}^{(m)}$ is the permutation matrix associated with the cyclic shift on $[m]$ by $+1$. Here $m+1=1$ in $[m]$.
\end{itemize}
Then, for every complex Banach space $E$ that is sphere equivalent to a complex Banach space of non-trivial type, we have that
\[
\mathcal{K}^{(\tau)}(\tilde{\mathbf{G}}_m,E)\succsim_{E,p} m^{-3}.
\]
\end{theorem}

The first very non-trivial point of this theorem is the left-hand side is \textit{strictly positive}; this follows from celebrated results of Lafforgue \cite{Lafforgue1}, \cite{Lafforgue2}, as we mentioned in Remark~\ref{remark=Lafforgue}. 
To estimate the decay in terms of $m$ from below, our argument is inspired by work of U. Hadad \cite{Hadad} employing bounded generation, originated in \cite{Shalom}. The following is the key to this part of argument.

\begin{lemma}[Dennis and Vaserstein, Lemma~$9$ in \cite{DennisVaserstein}] 
\label{lemma=DennisVaserstein}
Let $m\in \mathbb{N}_{\geq 3}$. Then every element in  $g\in \mathrm{SL}^{\pm}(m,\mathbb{F}_p[t])$ may be wriiten as 
\[
g=\gamma_1\gamma_1' \gamma_2\gamma_2' g_0,
\]
where $\gamma_1,\gamma_2$ belong to the group of all upper triangular matrices with $1$ on diagonal, $\gamma_1',\gamma_2'$ belong to the group of all lower triangular matrices with $1$ on diagonal, and $g_0$ lies in the copy of $\mathrm{SL}^{\pm}(3,\mathbb{F}_p[t])$ on the upper left corner.
\end{lemma}
This is because $\mathbb{F}_p[t]$ is a euclidean domain and in particular is Dedekind: For a Dedekind domain, the \textit{stable range} in the sense of Bass is at most $2$.

In the bounded generation argument, the following ``\textit{triangle inequality}'' plays an important role: For an isometric linear representation $\pi\colon G\to O(E)$, for every $\xi\in E$, it holds that for all $g,h\in G$,
\[
\|\pi(gh)\xi-\xi\|\leq \|\pi(gh)\xi-\pi(g)\xi\|+\|\pi(g)\xi-\xi\|=\|\pi(g)\xi-\xi\|+\|\pi(h)\xi-\xi\|.
\]

\begin{proof}[Proof of Theorem~$\ref{theorem=Tau}$]
Fix $m\in \mathbb{N}_{\geq 3}$. Let $S_1^{(m)}$ be the set of all elementary matrices of the form $e_{i,j}^a$ for $i\ne j \in [m]$ and $a\in \{\pm 1,\pm t\}$. Let $S_2^{(m)}$ be the set of all diagonal matrices in $\tilde{G}_m$ all of whose diagonal entries but one are $1$. Let $S^{(m)}=S_1^{(m)}\cup S_2^{(m)}$. Then, $S^{(m)}$ is a finite generating set of $\tilde{G}_m$. We also consider an \textit{infinite} subset $E^{(m)}$ of all elementary matrices in $\tilde{G}_m$.

Lafforgue's deep work in \cite{Lafforgue1}, \cite{Lafforgue2} implies the following. A more precise argument goes as follows: For a complex Banach space $F$ of non-trivial type, it follows from Lafforgue's work that $\mathrm{SL}(3,\mathbb{F}_p((t^{-1})))$ has strong Banach property $(\mathrm{T})$ for a certain class of Banach spaces that contains $L_2([0,1],F)$. It implies the fixed point property with respect to all affine isometric actions on $L_2([0,1],F)$. By $L_2$-induction (see for instance \cite[Section~8]{BFGM}), it then follows that an ($L_2$-integrable) lattice $\mathrm{SL}(3,\mathbb{F}_p[t])$ in $\mathrm{SL}(3,\mathbb{F}_p((t^{-1})))$ has the fixed point property with respect to ones on $F$. Since $\mathrm{SL}^{\pm}(3,\mathbb{F}_p[t])$ is an overgroup of $\mathrm{SL}(3,\mathbb{F}_p[t])$ of index $2$, so does $\mathrm{SL}^{\pm}(3,\mathbb{F}_p[t])$. By \cite[Theorem~1.3.(1)]{BFGM}, $\mathrm{SL}^{\pm}(3,\mathbb{F}_p[t])$ hence has property $(\mathrm{T}_F)$ in the sense of Bader--Furman--Gelander--Monod \cite{BFGM}. By recalling Remark~\ref{remark=ell_2}, we then observe that a certain ``Kazhdan constant'' associated to property $(\mathrm{T}_F)$ above is strictly positive. Thus, we obtain property $(\tau)$ for $\mathrm{SL}^{\pm}(3,\mathbb{F}_p[t])$ with respect to $F$. 
\begin{theorem}[V. Lafforgue]\label{theorem=Lafforgue}
The group $G=\mathrm{SL}^{\pm}(3,\mathbb{F}_p[t])$ has property $(\tau)$ with respect to every complex Banach space $F$ of non-trivial type. It means that, for every $F$ as above, it holds that
\[
\mathcal{K}^{(\tau)}(\mathbf{G},F)>0
\]
for all $($equivalently, some$)$ markings of $G$.
\end{theorem}

In \cite[Proposition~4.2]{MimuraSphereEquivalence}, the first named author showed that for a fixed (marked) group, the strict positivity of the $(\tau)$-type constant is preserved under replacing the Banach space with one that is sphere equivalent to the original one; see also Lemma~\ref{lemma=FiniteGroupQuotient} and discussions above it. Therefore, Theorem~\ref{theorem=Lafforgue} remains true if $F$ is replaced with $E$ as in Theorem~\ref{theorem=Tau}. Fix such an $E$ and 
define 
\[
\varepsilon_{E,p} =\mathcal{K}^{(\tau)}(\mathrm{SL}^{\pm}(3,\mathbb{F}_p[t]);S^{(3)}) >0.
\]
We first claim the following.
\begin{lemma}\label{lemma=claim1}
In the setting above, 
\[
\mathcal{K}^{(\tau)}((\tilde{G}_m;S^{(m)}),E)\geq \frac{1}{2m^2}\varepsilon_{E,p}.
\]
\end{lemma}

\begin{proof}[Proof of Lemma~$\ref{lemma=claim1}$]
We fix an isometric linear representation $\pi\colon \tilde{G}_m\to O(E)$ that factors through a finite group quotient. Take an arbitrary $\xi\in E$, and fix it.

Observe that there exists plenty of injective homomorphisms $\iota\colon \tilde{G}_3\hookrightarrow \tilde{G}_m$ that satisfies $\iota(S^{(3)})\subseteq S^{(m)}$; for instance, the natural embedding into the upper left corner satisfies it. For each such an $\iota$, we have that 
\[
\sup_{s\in S^{(m)}}\|\pi(s)\xi-\xi\|\geq \sup_{s\in \iota(S^{(3)})}\|\pi(s)\xi-\xi\|\geq \varepsilon_{E,p} \sup_{h\in \iota(\tilde{G}_3)} \|\pi(h)\xi-\xi\|.
\]
(Note that the restriction of $\pi$ on $\iota(\tilde{G}_3)$ also \textit{factors through a finite group quotient.}) We may find such $\iota$'s such that $\iota(\tilde{G}_3)$ for various $\iota$ altogether cover $E^{(m)}$. This implies that
\[
\sup_{s\in S^{(m)}}\|\pi(s)\xi-\xi\|\geq \varepsilon_{E,p} \sup_{h\in E^{(m)}} \|\pi(h)\xi-\xi\|
\]
and 
\[
\sup_{s\in S^{(m)}}\|\pi(s)\xi-\xi\|\geq \varepsilon_{E,p} \sup_{g_0\in \mathrm{SL}^{\pm}(3,\mathbb{F}_p[t])} \|\pi(g_0)\xi-\xi\|,
\]
where $\mathrm{SL}^{\pm}(3,\mathbb{F}_p[t])$ sits on the upper left corner of $\tilde{G}_m$.

We now recall Lemma~\ref{lemma=DennisVaserstein}. Observe that every $\gamma_1,\gamma_2,\gamma_1',\gamma_2'$ may be written as at most $(m(m-1))/2$ elements in $E^{(m)}$. Therefore, by the triangle inequality mentioned above Theorem~\ref{theorem=Lafforgue}, we conclude that
\[
\sup_{g\in \tilde{G}_m}\|\pi(g)\xi-\xi\|\leq 2m(m-1) \sup_{h\in E^{(m)}}\|\pi(h)\xi-\xi\|+\sup_{g_0\in \mathrm{SL}^{\pm}(3,\mathbb{F}_p[t])} \|\pi(g_0)\xi-\xi\|.
\]
By combining the three inequalities above, we obtain that
\[
\sup_{s\in S^{(m)}}\|\pi(s)\xi-\xi\|\geq \frac{\varepsilon_{E,p}}{2m(m-1)+1} \sup_{g\in \tilde{G}_m} \|\pi(g)\xi-\xi\|.
\]
This proves our claim.
\end{proof}

Finally, we claim that every element $s\in S^{(m)}$ may be written as at most $5000m$ products of elements in $\{\tilde{\sigma}^{(m)},\tilde{\upsilon}^{(m)},\tilde{\delta}^{(m)},\tilde{\tau}^{(m)}\}$ and their inverses. To show this, we employ a result of Kassabov and T. R. Riley \cite[Proposition~2.1]{KassabovRiley}, which implies that every element in $\tilde{G}_m$ of the form $e_{i,j}^a$, $i\ne j\in [m]$ and $a\in \{\pm 1\}$, may be written as at most $100m$ products of elements in $\{\tilde{\sigma}^{(m)},\tilde{\upsilon}^{(m)},\tilde{\delta}^{(m)},\tilde{\tau}^{(m)}\}$ and their inverses. Observe that for every $1\ne j\in [m]$, permutation matrix associated with the transposition on $\{1,j\}$ as the product of $3$ of such elements and $\tilde{\delta}^{(m)}$. From them, we may construct elements in $S_2^{(m)}$ and elements of the form $e_{i,j}^{\pm t}$, $i\ne j\in [m]$, within small numbers of products of elements above; compare with $(\#)$ in Subsection~\ref{subsection=ExtendedLinearGroups}. Hence, in the setting of the proof of Lemma~\ref{lemma=claim1}, we have that
\[
\sup_{s\in S^{(m)}}\|\pi(s)\xi-\xi\|\leq 5000m \cdot (\sup \{\|\pi(u)\xi-\xi\|:u\in \{\tilde{\sigma}^{(m)},\tilde{\upsilon}^{(m)},\tilde{\delta}^{(m)},\tilde{\tau}^{(m)}\}\}).
\]
Therefore,
\[
\mathcal{K}^{(\tau)}(\tilde{\mathbf{G}}_m,E)\geq  \frac{1}{10000m^3}\varepsilon_{E,p},
\]
as desired.
\end{proof}

\begin{proof}[Proof of Theorem~$\ref{theorem=BanachSpectralGap}$]
Combine Theorem~\ref{theorem=Tau} with Lemma~\ref{lemma=SpectralGapTau} (see also Remark~\ref{remark=ell_2}). Observe that for each $m\in \mathbb{N}_{\geq 3}$, $(G_m;\sigma^{(m)}, \upsilon^{(m)},\delta^{(m)},\tau^{(m)})$ is a finite marked group quotient of $(\tilde{G}_m;\tilde{\sigma}^{(m)},\tilde{\upsilon}^{(m)},\tilde{\delta}^{(m)},\tilde{\tau}^{(m)})$.
\end{proof}

\begin{proof}[Proof of Theorem~$\ref{mthm=PoorCompression}$]
Let $\Gamma_m=\mathrm{Cay}(G_m;\sigma^{(m)}, \upsilon^{(m)},\delta^{(m)},\tau^{(m)})$ and $E$ be as in Theorem~$\ref{mthm=PoorCompression}$. Since this graph is always $8$-regular and $\sharp(G_m)\succsim p^{n_m(m^2-m)}$, we have that
\[
\mathrm{diam}(\Gamma_m)\succsim m^2 \log n_m.
\]
For a given $\rho$, if necessary, replace it with a smaller proper function that satisfies the condition of Lemma~\ref{lemma=Austin} on $\rho$. Therefore, Theorem~\ref{theorem=JolissaintValette}, Lemma~\ref{lemma=Austin}, and Theorem~\ref{theorem=BanachSpectralGap} will provide an explicit (order of) $(n_m)_{m\in \mathbb{N}_{\geq 3}}$ such that $(\rho,\omega(r)=r) \not \in \mathcal{CP}_E(\bigsqcup_m\Gamma_m)$. 
For instance, if $\rho(r)\asymp (C'\log \log r) \vee 0 $ for some $C'>0$, then for every $C>0$, $e^{e^{e^{Cm^3}}}\precnsim n_m$  does the job.
\end{proof}

\begin{remark}\label{remark=BetaLessThan1/2}
Owing to the points that we discussed in Remark~\ref{remark=deLaatMimuradelaSalle}, we may obtain the same conclusion for $V=V_{(l_m)_m}$ as in Proposition~\ref{proposition=ChoiceOfGenerators} if the Banach space $E$ is sphere equivalent to a complex Banach space $F$ that satisfies ``$(1.1)$ for some $\beta<1/2$'' as in \cite{deLaatMimuradelaSalle}. The class of all such $F$ is included in the class of all complex Banach spaces of non-trivial type. It is not known whether this inclusion is strict; see also \cite[Example~4.11.$(4)$ and $(5)$]{MSPartII}.
\end{remark}

We finally remark that similar poor compression results were previously known by \cite{ArzhantsevaDrutusapir}, \cite{Austin} and \cite{OlshanskiiOsin} even in the setting of groups. Our point here is that our examples (graphs in $V$ and $Y$) appeared in a very natural way in a context of the question of \textit{unbounded rank expanders}; see \cite{LubotzkyWeiss} and \cite{BookLubotzky}. What we did here is not only to construct an example that admits coarse embeddings but only with poor compression functions; we discovered that certain natural examples, that have been paid attention by various researchers before consideration of the authors, have such property. We furthermore note that our examples are rather concrete: The example $V=V_{(l_m)_m}$ as in Proposition~\ref{proposition=ChoiceOfGenerators} is completely explicit. Our example $Y=Y_{p,(n_m)_m}$ is ``semi-explicit'' because in general to obtain a generator of $\mathbb{F}_{p^{n_m}}^{\times}$ explicitly is an issue; however, if we replace the coefficient ring from $\mathbb{F}_{p^{n_m}}$ to, say, $\mathbb{F}_p[t]/(t^{n'_m}-t)$ and if we consider $\mathbb{F}_p[t]\twoheadrightarrow \mathbb{F}_p[t]/(t^{n'_m}-t)$; $t\mapsto t$ instead of $\mathbb{F}_p[t]\twoheadrightarrow \mathbb{F}_{p^{n_m}}$; $t\mapsto t_{n_m}$, then this modification gives a totally explicit example $\overline{Y}=\overline{Y}_{p,(n_m')_m}$. This  $\overline{Y}$ remains to satisfy the conclusions of $(ii)$ of Corollary~\ref{mcor=ChoiceOfGenerators} and of Theorem~\ref{mthm=PoorCompression}; recall Remark~\ref{remark=OtherRings}.

\section*{Acknowledgments}
The authors wish to express their gratitude to Professor Guoliang Yu and Professor Qin Wang for their kind invitation to Fudan University in Shanghai in July, 2013. Part of this work was done during that stay. Some other part of this work was done during the  two-year stay of the first-named author in the \'{E}cole Polytechnique F\'{e}d\'{e}rale de Lausanne supported by Grant-in-Aid for JSPS Oversea Research Fellowships. The first-named author wishes to express his gratitude to Professor Nicolas Monod and Mrs. Marcia Gouffon for their hospitality and help on his visit. The two authors  are grateful to Romain Tessera on LEF groups and for drawing their attention to the example as in Lemma~\ref{lemma=Tessera}, and Ana Khukhro and Yves Stalder on wreath products. They also thank Narutaka Ozawa for discussions, and Yuhei Suzuki for the careful reading of the first draft of this paper. The first named author is grateful to Rostislav I. Grigorchuk for providing him with Example~\ref{example=Grigorchuk}, to Goulnara Arzhantseva for discussions on alternating groups and suggestion for adding Remark~\ref{remark=OtherRings}, to Laurent Bartholdi and Anna Erschler for discussions on \cite{BartholdiErschler}, to Yves de Cornulier and Martin Kassabov respectively for the reference \cite{BookMalcev} and  \cite{Stepin}, and \cite{LubotzkyWeiss}, and to Vsevolod Gubarev for the terminology ``unimodular'' for $\mathrm{SL}^{\pm}$-group. The authors are grateful to the referee for several comments and for supplying them the references  \cite{DelabieKhukhro} and \cite{Guentner}.

\bibliographystyle{amsalpha}
\bibliography{mimura_sako_cayley1.bib}

\end{document}